\documentclass[review]{elsarticle}

\journal{Journal of \LaTeX\ Templates}

\usepackage{amssymb}
\usepackage{latexsym}
\usepackage{mathtools}
\usepackage{url}
\usepackage{xcolor}

\usepackage{amsmath}
\usepackage{mathtools}
\usepackage{amsfonts}
\usepackage{amsthm}
\usepackage{bm}
\usepackage{physics}

\newtheorem{theorem}{Theorem}

\newtheorem{lemma} {Lemma}

\theoremstyle{remark}

\theoremstyle{definition}
\newtheorem{definition}{Definition}[section]

\usepackage{floatrow}

\usepackage[capitalise]{cleveref}
\graphicspath{ {images/} }

\usepackage{algpseudocode}
\usepackage{algorithm}
\usepackage[capitalise]{cleveref}
\crefformat{equation}{(#2#1#3)}
\crefname{line}{Algorithm}{algorithm}
\crefname{Item}{}{}
\usepackage{tikz-cd}

\newcommand{\mat}[1]{\mathbb{#1}}

\newcommand{\dual}[1]{\tilde{#1}}

\newcommand{\R}{\mathbb R}

\newcommand{\mc}{\mathcal}

\newcommand{\inner}[2]{\langle {#1} \, , {#2} \rangle}
\newcommand\chain[1]{C_{#1}(\mc K)}
\newcommand\cochain[1]{C^{#1}(\mc K)}

\newcommand\bd[1]{\partial_{#1}}
\newcommand\cobd[1]{\delta^{#1}}

\newcommand\bbd[2]{\mathrm{bd}_{#1}{(#2)}}
\newcommand\bcobd[2]{\mathrm{cobd}_{#1}{(#2)}}

\newcommand\block[2]{{{#1}|_{#2}}}

\newcommand\vc{\mathfrak{v}}
\newcommand\ec{\mathfrak{e}}
\newcommand\fc{\mathfrak{f}}
\newcommand\cc{\mathfrak{c}}

\newcommand\Cbasis{\widehat{\mc B}}
\newcommand\basis{\mc B}
\newcommand\degree[1]{\mathrm{deg}(#1)}
\newcommand\ob{(C^*,\mc B)}
\algdef{SE}[DOWHILE]{Do}{doWhile}{\algorithmicdo}[1]{\algorithmicwhile\ #1}%

\newtheorem{example}[theorem]{Example}

\begin{document}

\begin{frontmatter}

\title{Inverting the discrete curl operator: a novel graph algorithm to find a vector potential of a given vector field}

\author[1]{Silvano Pitassi\corref{cor1}}
\ead{pitassi.silvano@spes.uniud.it}
\cortext[cor1]{Corresponding author:
  Tel.: +039-0432-558037;}
\author[2]{Riccardo Ghiloni}
\ead{riccardo.ghiloni@unitn.it}
\author[1]{Ruben Specogna}
\ead{ruben.specogna@uniud.it}

\address[1]{University of Udine, Polytechnic Department of Engineering and Architecture, EMCLab, via delle scienze 206, 33100 Udine, Italy}
\address[2]{University of Trento, Department of Mathematics, via Sommarive 14, 38123 Povo-Trento, Italy}

\begin{abstract}
We provide a novel framework to compute a discrete vector potential of a given discrete vector field on arbitrary polyhedral meshes.
The framework exploits the concept of acyclic matching, a combinatorial tool at the core of discrete Morse theory.
We introduce the new concept of \emph{complete acyclic matchings} and we show that they give the same end result of Gaussian elimination.
Basically, instead of doing costly row and column operations on a sparse matrix, we compute equivalent cheap combinatorial operations that preserve the underlying sparsity structure.
Currently, the most efficient algorithms proposed in literature to find discrete vector potentials make use of tree-cotree techniques.
We show that they compute a special type of complete acyclic matchings.
Moreover, we show that the problem of computing them is equivalent to the problem of deciding whether a given mesh has a topological property called \emph{collapsibility}.
This fact gives a topological characterization of well-known termination problems of tree-cotree techniques.
We propose a new recursive algorithm to compute discrete vector potentials.
It works directly on basis elements of $1$- and $2$-chains by performing elementary Gaussian operations on them associated with acyclic matchings.
However, the main novelty is that it can be applied recursively.
Indeed, the recursion process allows us to sidetrack termination problems of the standard tree-cotree techniques. 
We tested the algorithm on pathological triangulations with known topological obstructions.
In all tested problems we observe linear computational complexity as a function of mesh size.
Moreover, the algorithm is purely graph-based so it is straightforward to implement and does not require specialized external procedures.
We believe that our framework could offer new perspectives to sparse matrix computations.
\end{abstract}

\end{frontmatter}

%\linenumbers

\section{Introduction}
In this paper we present a novel framework to solve the discrete version of the following potential problem: determine a vector field with specified curl, i.e. a vector potential of a given vector field.
To do so, an ``inverse curl'' operator is required, necessarily including the specification of additional conditions (a gauge condition in physical parlance) to uniquely define a resulting solution.
Let us now formally introduce the problem by first considering the continuous case and then moving towards its discrete version, where we define discrete counterparts of vector fields and differential operators acting on them.

Let us consider a bounded domain $\Omega$ of $\R^3$, where a vector field $\bm J$ is defined.
Assume that $\Omega$ is topologically trivial, i.e., it is homeomorphic to a closed $3$-ball (or, equivalently, to a cube).
We first determine necessary and sufficient conditions for assuring that $\bm J$ is the curl of a vector field $\bm H$, i.e. $\bm H$ is a vector potential of $\bm J$.
The answer is well-known, being a classical problem in vector analysis.
A vector field $\bm J$ is the curl of a vector field if and only if its divergence is zero and its flux is vanishing across all the (but one) connected components of $\partial\Omega$.
In our case, $\partial \Omega$ is connected, so the flux condition can be omitted since it is automatically verified.

Let us consider a \emph{mimetic discretization} of our continuous problem.
Mimetic discretization methods like the Mimetic Finite Difference method (MFD) \cite{mimeticjcp}, Discrete Geometric Approach (DGA) \cite{pitassijcp2021}, Finite Integration Technique (FIT), Discrete de Rham (DDR) methods \cite{dipietrojcp2021} include the structure of exterior calculus, thus retaining fundamental properties of the continuous theory.

We recall the main concepts of mimetic discretization in \cref{chains}.
We cover the domain $\Omega$ of $\R^3$ with a polyhedral mesh $\mc K$, namely a regular CW cell complex having cells being polyhedra.
Continuous vector fields and differential operators are replaced by their discrete counterparts.
We define a discrete vector field as a collection of \emph{degrees of freedom} (DoFs).
We introduce the discrete vector field $\bm h$, an array of DoFs with values on each edge of $\mc K$; next, the discrete vector field $\bm i$, an array of DoFs with values on each face of $\mc K$.
The central problem of this paper is to find a discrete vector potential $\bm h$ such that
\begin{equation}
\mat C \bm h = \bm i,
\label{main1}
\end{equation}
where $\mat C$ is the usual incidence matrix between the faces and the edges of the mesh.
In intimate analogy with the continuous case, a necessary and sufficient condition to have a consistent discrete vector potential $\bm h$ is that array $\bm i$ represents a \emph{discrete solenoidal} vector field, i.e. it verifies $\mat D \bm i = \bm 0$, where $\mat D$ is the incidence matrix between cells and faces.

Our main motivation to solve problem \cref{main1} stems from the fact that this algorithmic primitive is an enabling technology for solving many problems arising in computational physics, from electromagnetism to elasticity and fluid mechanics \cite{jmaa1990}. First, it can be used to solve the vector laplacian in nearly linear time \cite{Cohen2014Solving1I}. The idea is that, instead formulating the vector laplacian by using a vector potential, the scalar potential \cite{Webb1989ASS} or the mixed-hybrid \cite{pitassijcp2021} formulations could be used instead, which produce linear systems that can be solved in nearly linear time by using algebraic multigrid methods. Second, inverse discrete curl is at the root of efficient algorithms to compute a cohomology basis and source fields for solving magnetostatics and eddy current problems by mimetic or finite element methods \cite{clenet1998,sinum1,cpcds,vallijcp2015}.
We think $\bm h$ as a discrete magnetic field and $\bm i$ as a discrete current; then \cref{main1} expresses the so-called discrete Amp\`{e}re's law.
In the electromagnetic literature, discrete fields $\bm h$ satisfying \cref{main1} are often called \emph{source fields}.
A different application in computational electromagnetics is to find a magnetic vector potential from a magnetic induction field \cite{jcp2016_adams}.

There are two analogous discrete potential problems: the problem of determining a scalar potential with assigned gradient and a vector field with assigned divergence.
However, as we will see in our discourse, and as been already pointed out in literature \cite{Rodrguez2015FiniteEP}, these two problems are less challenging than problem \cref{main1}, since they are easily solved in linear worst-case complexity using standard spanning tree constructions.

Many different algorithms have been proposed in literature to solve linear system \cref{main1}.

A naive solution to the problem would be to solve \cref{main1} by a linear system solver.
However, these techniques are not feasible from the practical point of view since they show cubic worst-case computational complexity.
Instead, the most efficient methods to solve linear system \cref{main1} are based on the so-called \emph{tree-cotree decomposition}.
Tree-cotree decomposition arises from graph theory and consists in partitioning the edges of a graph into a \emph{spanning tree} and its complement, referred to as the cotree.
The basic idea, rooted in the works \cite{Webb1989ASS, clenet1998}, goes as follows.
First, the values of array $\bm h$ relative to a \emph{spanning tree} on the vertices and edges of $\mc K$ are set to zero.
Next, discrete Amp\`{e}re's law is  iteratively enforced on each face of the mesh.
This principle is the core of the algorithms proposed in \cite{Dlotko2010CriticalAO, Rodrguez2015FiniteEP}.

Tree-cotree techniques are frequently claimed to be general in the literature.
However, proofs about their termination are not discussed at all \cite{Webb1989ASS,clenet1998}.
In fact, such techniques are not guaranteed to converge.
A careful analysis of the termination properties can be found in \cite{Dlotko2010CriticalAO}, where it is shown that they strongly depend on the choice of the spanning tree.
In particular, there exist situations involving a topologically trivial complex (for example, certain meshes of a cube or a 3-ball) and a spanning tree on which such techniques do not terminate.

To overcome termination issues of tree-cotree techniques, different approaches have been proposed.
The approach in \cite{Dlotko2011EfficientGS} is based on the idea of symbolic computations.
Although being general, this approach is slower than previous approaches.
Moreover, it is difficult to implement in practice since it requires specialized data structures to manage the symbolic computations; for instance, it requires object oriented programming languages to implement the symbolic computations.
In \cite{Rodrguez2015FiniteEP}, termination issues are solved using an explicit formula based on a double integral computation.
However, these double integral computations require specialized algorithms that are time consuming.

The aim of this paper is to describe an efficient and easy to implement algorithm to solve linear system \cref{main1}.
Our novel algorithm is based on concepts of \emph{discrete Morse theory} \cite{Forman1998MorseTF}.
This theory employs a construction called \emph{acyclic matching} \citep{Kozlov2008}, which collects combinatorial operations analogous to topological operations in the continuum.
In \Cref{discrete_morse} we review basic concepts at the core of discrete Morse theory as well as our specialized definitions.

Starting from acyclic matchings of discrete Morse theory, we derive three main contributions.

Firstly, we provide a unified framework based on discrete Morse theory for the solution of linear system \cref{main1} and we present in \cref{gaussian}.
We introduce the novel concept of \emph{complete acyclic matching}.
The crucial fact of our framework is the following: complete acyclic matchings give the same end result of Gaussian elimination.
What the whole procedure boils down to is avoiding costly matrix algebra operations by performing equivalent cheap combinatorial operations.
By this procedure, the sparsity of linear system \cref{main1} is retained, whereas this would not be
the case in the standard Gaussian elimination.
In fact, sparse systems become dense in intermediate steps due to an inconvenient choice of the backward elimination algorithm.

Secondly, using the newly introduced framework, we show that tree-cotree decomposition techniques are algorithms to compute specialized complete acyclic Morse matchings.
Although they use a different language to describe the same actions, they provide the same end results.
We show that the question of finding this kind of specialized complete acyclic matchings is equivalent to the following \emph{collapsibility problem}: decide whether a 3-dimensional simplicial complex embedded in $\R^3$ with trivial topology and with Lipschitz boundary is \emph{collapsible}.
This fact establish the topological nature of termination problems of tree-cotree techniques and shows why they arise in practical applications.
Indeed, there exist examples of triangulations of 3-balls that are not collapsible \cite{Lutz2013KnotsIC}.
As a result, for certain triangulations, avoiding termination problems by a careful choice of the input spanning tree of tree-cotree techniques is an impossible task and thus we need to resort to new approximation strategies.

Thirdly, we provide a new recursive algorithm to compute discrete vector potentials we describe in \cref{algorithm}.
It consists of a new greedy heuristics to construct acyclic matchings together with a recursive construction, which has no analogous in classical discrete Morse theory.
Indeed, we do not define a new chain complex, the so-called \emph{Morse complex}, but instead we employ acyclic matchings to express new basis elements during Gaussian elimination in terms of the original basis elements of the vector space of 1-chains and the vector space of 2-chains.
Fundamentally, we do not tackle the problem of computing specialized complete acyclic matchings like in tree-cotree techniques, which, as discussed above, suffers from well-known topological termination problems.
Instead, we introduce a recursive approach whose basic outline is as follows. 
We first try to find a complete acyclic matching on linear system \cref{main1}.
If we do not succeed, we transform linear system \cref{main1} into a new smaller linear system, by considering a suitable subset of the newly computed basis during Gaussian elimination. 
We show that if we can solve this new smaller linear system, then we can get a solution of linear system in \cref{main1}.
The crucial fact is that we can recursively apply the algorithm on the new smaller linear system.
More specifically, we try to find a complete acyclic matching on the corresponding new smaller linear system.
Again, acyclic matchings provides a new basis during Gaussian elimination so that the whole approach can be applied recursively.
The fundamental computational advantage of our recursive approach is that all the costly algebraic operations are replaced by elementary cheap combinatorial operations on significantly smaller instances of the original linear system \cref{main1}.

We tested our recursive algorithm on challenging benchmark problems.
These include also pathological meshes with known topological obstructions.
Our algorithm exhibits linear computational complexity for all tested problems.
Moreover, the algorithm is also purely graph-theoretic, so straightforward to implement and do not require any additional specialized procedure.
We collect simulation results in \cref{numerical}.

While obtaining a theoretical linear worst-case complexity bound is hard, our algorithm solves all issues that are typically found in practice.
Moreover, the whole framework proves to be general and offers new perspectives to sparse matrix computations.
These observations are summarized in \cref{conclusions}.

\section{Notation}
\label{chains}
The domain of interest of this paper is a closed and bounded polyhedral domain $\Omega$ of $\R^3$ with Lipschitz boundary. We assume that $\Omega$ has trivial topology, i.e., it is homeomorphic to a closed $3$-dimensional ball or, equivalently, $\Omega$ is simply connected and its boundary $\partial\Omega$ is connected ($\partial\Omega$ is homeomorphic to a $2$-sphere indeed); see \cite{Cantarella2002VectorCA} (Section 6) and \cite{Benedetti2010TheTO} (Section 3).
We consider a \emph{polyhedral cell complex} (or \emph{polyhedral mesh})  subdivision $\mathcal{K}$ of $\Omega$.
Elements of $\mc K$ are called \emph{cells}.
A $k$-cell $\sigma$ is a $k$-dimensional subset in $\R^3$ homeomorphic to a closed $k$-dimensional ball.
A $0$-cell is a point of $\R^3$.
We equip each $k$-cell with an inner orientation.
We denote by $\dim {\sigma}$ the dimension of the cell $\sigma \in \mc K$ and we write $\mc K_{k}$ the subcollection of all $k$-cells in $\mc K$.
We focus on the 3-dimensional case, thus we have 3-cells (or volumes) in $\mc K_3$, 2-cells (or faces) in $\mc K_2$, 1-cells (or edges) in $\mc K_1$ and 0-cells (or vertices) in $\mc K_0$.
We denote by $c$ a generic volume, by $f$ a face, by $e$ an edge and by $v$ a vertex.
We denote by $\vc, \ec, \fc, \cc$ the cardinality of $\mc K_{0},\mc K_{1}, \mc K_{2}, \mc K_{3}$, respectively.
A polyhedral cell complex $\mc K$ is \emph{simplicial} if all its cells are simplicies and the boundary of each cell has the natural simplicial decomposition, see \cite{Christiansen2008ACO}.
If $\mc K$ is simplicial, $\mc K$ is also called a triangulation of $\Omega$.

The mesh $\mc K$ has the structure of a \emph{(regular) cell complex}, namely, the following three conditions hold \cite{Christiansen2008ACO}.
First, for each $k$-cell $\sigma$ in $\mc K$ its \emph{boundary} $\partial \sigma$ is a union of $(k-1)$-cells in $\mc K$ for $k\in\{1,2,3\}$ \cite{Christiansen2008ACO}.
Second, given distinct $k$-cells $\sigma,\tau$, their intersection $\sigma \cap \tau$ is either empty or is a union of lower dimensional cells in $\mc K$.
Third, given a $l$-cell $\sigma$ and $k$-cell $\tau$ with $l \leq k$, $\sigma \neq \tau$ and $\sigma \cap \tau \neq \emptyset$, we have $\sigma \cap \tau \subset \partial \tau$.

We can now define a new object, called a real $k$-chain.
A $k$-chain $a$ of $\mathcal{K}$ is a formal linear combination of $k$-cells $a = \sum_{i=1}^r a_i \sigma_i$, where $\sigma_i$ are $k$-cells in $\mc K$ and $a_i$ are real coefficients.
The number $r$ denotes the cardinality of the collection of $k$-cells in $\mathcal{K}$ and is any number among $\vc, \ec, \fc$ or $\cc$.
The set of $k$-chains, equipped with the natural addition and scalar multiplication, provides a real vector space.
We denote it by $\chain{k}$.
Note that each $k$-cell is also a $k$-chain.
If $\sigma$ is a $k$-cell, by $-\sigma$ we denote the cell $\sigma$ but with opposite orientation.
The set of all $k$-cells form a basis for $\chain{k}$, which we call \emph{canonical basis} for $\chain{k}$.
We identify the boundary of each $k$-cell with the linear combination of the $(k-1)$-cells in $\partial \sigma$ defined by setting
\begin{equation}
\partial \sigma \coloneqq \sum_{i=1}^r w_i\,\rho_i,
\label{def_bd}
\end{equation}
where $w_i$ is different from zero if and only if $\rho_i \subset \partial \sigma$ and in this case, $w_i$ is equal to $+1$ if $\rho_i$ has the orientation induced by that of $\sigma$ by using the right-hand rule and $-1$ otherwise \cite{Tonti2013TheMS}.

The real vector space of $k$-chains and the real vector space of $(k-1)$-chains are connected by a linear map called \emph{boundary operator} $\bd{k}: \chain{k} \to \chain{k-1}$.
We define the boundary operator by linearity on the space of chains by setting
\begin{equation}
\bd{k} a \coloneqq \sum_{i=1}^r a_i \, \partial \sigma_i,
\label{bd_def}
\end{equation}
where $a = \sum_{i=1}^r a_i \sigma_i$ as above.
Note that \cref{bd_def} is well-defined since $\mc K$ is a cell complex.
Since $\mc K$ is a cell complex, it can verified that $\bd{k-1}\circ \bd{k}=0$ for $k\in\{1,2,3\}$; see, for example \cite{Christiansen2008ACO}.

Let us now consider the concept of a \emph{$k$-cochain}.
A $k$-cochain $b$ acts on a $k$-chain to produce a real number and therefore $k$-cochains are elements of the dual space of $\chain{k}$.
We define the vector space of $k$-cochains $\cochain{k}$ to be the dual space of linear functionals $b:\chain{k}\to \R$.
We denote the value of a $k$-chain $a$ under a $k$-cochain $b$ as $\inner{b}{a} \coloneqq b(a)$.
Let us consider the canonical basis $\{\sigma_i \in \chain{k} \mid i=1,\dots,r\}$ of the vector space of $k$-chains $\chain{k}$.
From basic linear algebra, there exist unique linear functionals $\{\sigma^i \in \cochain{k} \mid i=1,\dots,r\}$ such that
\begin{equation}
\inner{\sigma^i}{\sigma_j}=\delta_{ij},
\label{def_delta}
\end{equation}
where $\delta_{ij}$ is the Kronecker delta.
The set $\{\sigma^i \in \cochain{k} \mid i=1,\dots,r\}$ defined by \cref{def_delta} form a basis for the vector space of $k$-cochains $\cochain{k}$, which is called \emph{canonical dual basis}.
We have established a one-to-one correspondence between chains and cochains.
This chain-cochain natural duality yields the real linear isomorphism $\phi_k:\chain{k} \to \cochain{k}$ sending each $\sigma_i$ to $\sigma^i$.
We will write a generic $k$-cochain $b \in \cochain{k}$ as a sum $b=\sum_{i=1}^r b_i \sigma^i$ with real coefficients $b_i$.

For $k$-cochains, in intimate analogy with chains, we can define a \emph{coboundary operator} $\delta^k: \cochain{k} \to \cochain{k+1}$ as the \emph{dual of the boundary operator}, i.e., it is defined by requiring that, for every $b \in \cochain{k}$ and $a \in \chain{k+1}$, the following identity holds
\begin{equation}
\inner{\cobd{k}b}{a}=\inner{b}{\bd{k+1}a}.
\label{def_cobd}
\end{equation}
In mimetic methods the coboundary operator $\cobd{k}$ acts as a discrete counterpart of the continuous differential operators \cite{Lipnikov2014MimeticFD,Tonti2013TheMS}.
Specifically, $\cobd{0}$ acts as the discrete gradient, $\cobd{1}$ as the discrete curl and $\cobd{2}$ as the discrete divergence.

A straightforward calculation using \cref{def_cobd} shows that $\cobd{k} \circ \cobd{k-1}=0$ for $k\in\{1,2\}$.
These relations mimic the structure of continuous differential operators \cite{Lipnikov2014MimeticFD,Tonti2013TheMS}.
In particular, discrete differential operators form a chain complex
\begin{equation}
C^* = \cdots \xrightarrow{\cobd{k-1}} \cochain{k} \xrightarrow{\cobd{k}} \cochain{k+1}
\xrightarrow{\cobd{k+1}} \cdots,
\label{chain_complex}
\end{equation}
where $\cochain{k}=0$ if $k<0$ or $k>3$.
Since the domain $\Omega$ is topologically trivial the sequence is \emph{exact} for $k \neq 0$, i.e. it satisfies $\mathrm{im}(\cobd{k-1})=\mathrm{ker}(\cobd{k})$ for $k\neq0$.

In the case of $k$-chains, there is a natural choice of a basis given by the canonical basis.
Using the isomorphism $\phi_k: \chain{k} \to \cochain{k}$, we have also fixed a canonical dual basis for $\cochain{k}$.
Since the coboundary operator $\cobd{k}$ is a linear map between $\cochain{k}$ and $\cochain{k+1}$, it can be represented, using the fixed bases of $\cochain{k}$ and $\cochain{k+1}$, as a matrix.
Thus, to represent the coboundary operator as a matrix, we must always explicitly state which bases are chosen and, in fact, we will soon see the benefits of
changing the bases.

Let us now consider an arbitrary basis for the vector space of $k$-chains $\chain{k}$.
We denote it by $\basis_k=\{\dots,\xi_i,\dots\}$.
Using the isomorphism $\phi_k:\chain{k} \to \cochain{k}$, $\phi_k(\basis_k)=\{\dots,\xi^i = \phi_k(\xi_i),\dots\}$ is a basis for $\cochain{k}$.
In what follows, we take this process of choosing a basis for $\cochain{k}$ for granted.
When this is done, we say that we have chosen a basis $\basis = \bigcup_k \basis_k$ for the entire chain complex $C^*$, i.e. a basis for each $\cochain{k}$.
We write $\ob$ to denote a chain complex with a basis.
We define the \emph{canonical basis} $\Cbasis$ for $C^*$ to be the basis $\Cbasis=\bigcup_k \Cbasis_k$ where each $\Cbasis_k$ is the canonical basis for $\chain{k}$.

Having chosen bases in $\chain{k}$ and hence in $\cochain{k}$, we denote by $\mat D_k$ the matrix associated with $\delta^k$ for $k\in \{0,1,2\}$.
We introduce a more common notation from mimetic methods $\mat G \coloneqq \mat D_0$, $\mat C \coloneqq  \mat D_1$ and $\mat D \coloneqq  \mat D_2$.
We can also define $\mat D_3$ as the null operator from $\chain{3}$ to $0$.

We represent $k$-chains and $k$-cochains by vectors of size $r$ that contain the real numbers with respect to the ordered bases.
A $k$-chain $a=\sum_{i=1}^r a_i \xi_i$ in a basis $\{\dots, \xi_i,\dots\}$ is represented by the column vector $\R^r \ni \bm a = (a_1 \cdots a_r)^T$ and a $k$-cochain $b=\sum_{i=1}^r b_i \xi^i$ in the basis $\{\dots, \xi^i = \phi_k(\xi_i),\dots\}$ is represented by the column vector $\R^r \ni \bm b =(b_1, \cdots b_r)^T$.

As explained in the introduction, the aim of this paper is to devise an efficient algorithm for the solution of the following problem: find array $\bm h \in \R^\ec$ such that
\begin{equation}
\mat C \bm h = \bm i,
\label{main}
\end{equation}
where $\bm i \in \R^\fc$ satisfy $\mat D \bm i = \bm 0$.

The rank of the matrix $\mat C$ is not maximal, thus linear system \cref{main} has an infinite number of solutions.
In fact, if $\bm h$ is a solution of \cref{main}, then array $\bm h' = \bm h + \mat G \bm \psi$ is also a solution of \cref{main} since we have
\begin{equation}
\mat C \bm h' = \mat C (\bm h + \mat G \bm \psi ) = \mat C \bm h = \bm i,
\end{equation}
where $\bm \psi \in \R^\vc$ and we have used the chain complex property $\mat C \mat G = 0$ of \cref{chain_complex}.

In what follows, the concept of a partition of a given index set will play an important role.
Let $I$ be a finite index set.
A partition of $I$ is a family of disjoint subsets $\{I_1,\dots,I_p\}$ of $I$ such that $\bigcup_{l=1}^p I_l = I$.
The \emph{subvector of $\bm v = (v_i)_{i \in I} \in \R^I$ induced by $I_l$} is
\begin{equation}
\block{\bm v}{I_l} \coloneqq (v_i)_{i \in I_l},
\label{block}
\end{equation}
for $l\in \{1,\dots,p\}$.
A representation of the vector $\bm v$ as a block vector is given by 
\begin{equation}
\bm v = (\block{\bm v}{I_l})_{l \in \{1,\dots,p\}}.
\end{equation}
Let us consider a product of index sets $I$ and $J$.
We will need a corresponding notion of \cref{block} for a matrix whose entries are indexed by elements in $I\times J$.
Let us consider partitions of $I$ and $J$ as $\{I_1, \dots, I_p\}$ and $\{J_1, \dots, J_q\}$, respectively.
We have a corresponding partition of $I \times J$ as a family of disjoint subsets $\{O_1, \dots,O_n\}$ such that $O_i=I_l \times J_m$ for some $l\in\{1,\dots,p\}$, $m\in\{1,\dots,q\}$ and $I \times J = \bigcup_{i=1}^n O_i$.
The \emph{submatrix of $\mat A = (\mat A_{i,j})_{i \in I,j \in J} \in \R^{I \times J}$ induced by $I_l \times J_m$} is
\begin{equation}
\block{\mat A}{I_l \times J_m} \coloneqq (\mat A_{i,j})_{i \in I_l, j \in J_m} \in \R^{I_l \times J_m}.
\end{equation}
A representation of the matrix $\mat A$ as a block matrix is given by
\begin{equation}
\mat A = (\block{\mat A}{I_l \times J_m})_{l \in \{1,\dots,p\}, m\in\{1,\dots,q\}}.
\end{equation}

\section{Discrete Morse Theory}
\label{discrete_morse}

The underlying principle of our construction follows an ad hoc reformulation of Forman's Discrete Morse theory \cite{Forman1998MorseTF} given by Kozlov \cite{Kozlov2008}, where the basic tool is a combinatorial object called acyclic matching.
Several special cases of our construction have already appeared in literature.
We present a formulation of discrete Morse theory adapted to our purposes, along with smaller, more illustrative instances, which will provide insights on the structure of our algorithm.

\subsection{Informal introduction to discrete Morse theory}
\label{informal}
\begin{figure}[t]
    \centering
    \includegraphics[scale=0.8]{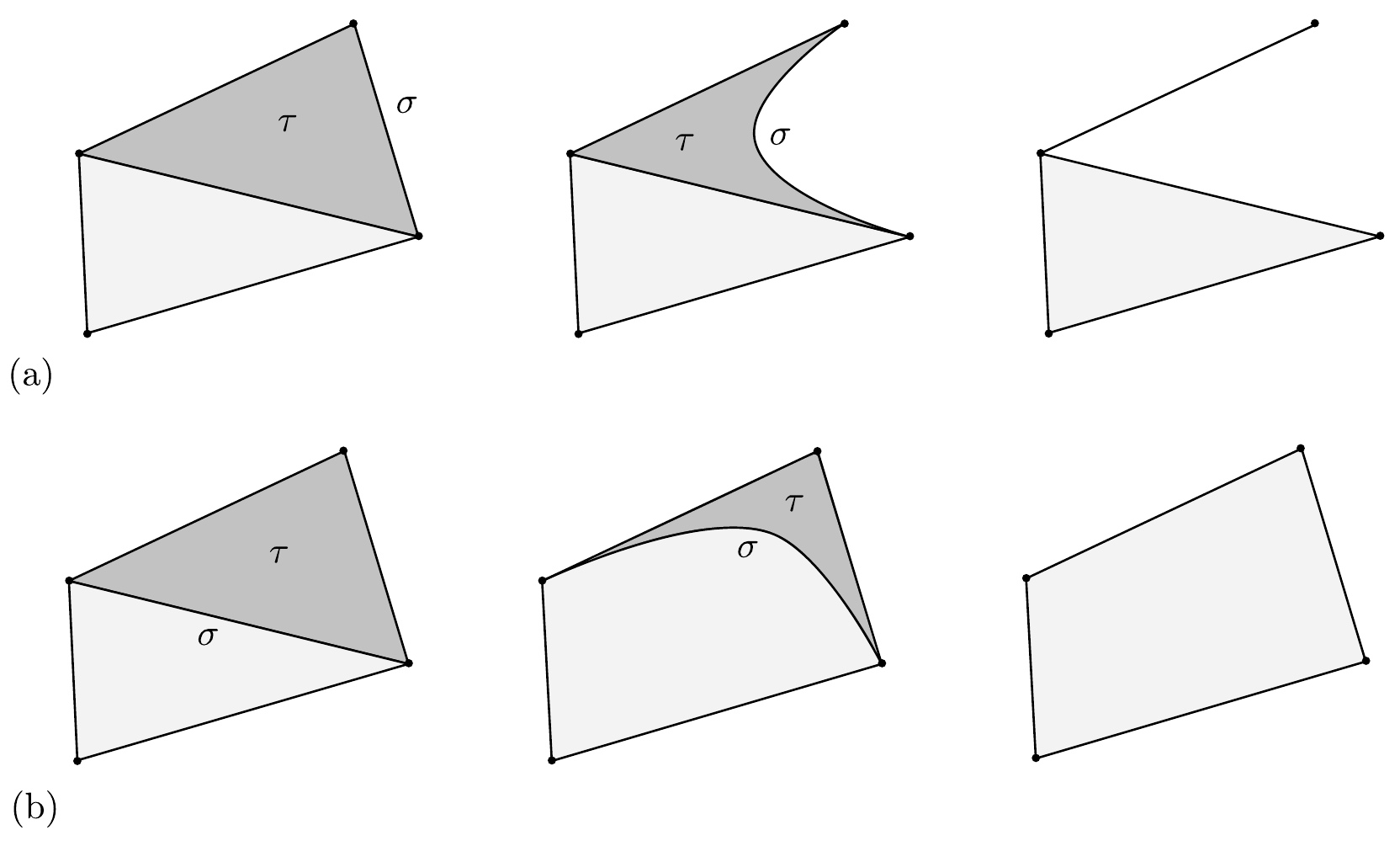}
    \caption{(a) Elementary collapse of free pair $(\sigma,\tau)$. (b) Internal collapse of pair $(\sigma,\tau)$; the resulting cell complex is not more simplicial.}
    \label{collapses}
\end{figure}

The first concept is that of \emph{elementary collapse}.
One may view discrete Morse theory as a generalization of the theory of simplicial collapses.
The concept of collapse, originated in Whitehead's work \cite{Whitehead1939SimplicialSN}, provides a combinatorial operation that is analogous to the continuous operation called deformation retraction, i.e., the operation of continuously shrinking a topological space to a subset.
More specifically, let $(\sigma,\tau)$ be a pair of cells such that $\sigma \subset \tau$ and $\dim \sigma = \dim \tau -1$.
For this pair, to induce an elementary collapse, we require $\tau$ to be a cell of maximal dimension in $\mc K$ and the only one cell of $\mc K$ containing $\sigma$; we refer to this as saying that the pair $(\sigma,\tau)$ is \emph{free} in $\mc K$.
Equivalently, we also say that $\sigma$ is free in $\mc K$; see \cref{collapses}(a). 
We say that $\mc K$ \emph{collapses to} $\mc L$ if one could get from $\mc K$ to $\mc L$ in a finite sequence of elementary collapses.
If $\mc K$ is equivalent to a single vertex, then we say that $\mc K$ is \emph{collapsible}; in this case there exists a sequence of elementary collapses leaving a single vertex.

Dropping the uniqueness condition on $\tau$, we obtain what we refer to as an \emph{internal collapse}, see \cref{collapses}(b).

Geometrically, in both cases, we obtain a collapse of the pair $(\sigma,\tau)$ by contracting the whole cell $\tau$ onto $\partial \tau \setminus \sigma$.

In intimate analogy with elementary collapses, we may combine many internal collapses to form a sequence of internal collapses, again without affecting the homotopy type.

We thus have a family of pairs $\{(\sigma_1,\tau_1),\dots,(\sigma_n,\tau_n)\}$ to be collapsed, in this order.
One may view the set of all such pairs as a \emph{matching} on $\mc K$.
Accordingly, we refer to cells contained in some pair as \emph{matched} and other cells as \emph{unmatched} or \emph{critical} (with respect to the matching).

Let $\mc K^{\,(i)}$ be the resulting cell complex after the first $i$ collapses.
For the pairs to form a sequence of elementary collapses, we require that each new pair $(\sigma_i,\tau_i)$ is free in $\mc K^{\,(i-1)}$.
For generic collapses we apply the same requirement, except that we restrict our attention to the family of matched cells.
Specifically, we do not require $(\sigma_i,\tau_i)$ to be free in $\mc K^{\,(i-1)}$, but $\tau_i$ must be the only matched cell of $\mc K^{\,(i-1)}$ containing $\sigma_i$.
Equivalently, for each $i$, we should have that $\sigma_i$ is not contained in $\tau_{i+1},\dots,\tau_n$ for $i \in \{1,\dots,n-1\}$.
We refer to a matching on $\mc K$ admitting an ordering with this property as \emph{acyclic}.
We formalize all this concepts in \cref{acyclic_matchings}.

The main theorem of discrete Morse theory states that an acyclic matching induces a homotopy equivalence between $\mc K$ and the so-called Morse complex, a cell complex formed by critical cells only \cite{Kozlov2008} (Theorem 11.13 (b)).

\subsection{Acyclic matchings}
\label{acyclic_matchings}
We start our exposition by examining acyclic matchings from a purely combinatorial point of view without any reference to topology.
Indeed, our interest is in using discrete Morse theory to develop a fast algorithm for the solution of linear system \cref{main} to be applied to cell complexes arising from experimental or numerical meshes of real case problems.
We give a specific version of combinatorial discrete Morse theory by Kozlov \cite{Kozlov2008} that is adapted to our purposes.

For any $\sigma, \tau \in \Cbasis_k$, define $\inner{\sigma}{\tau}$ to be $1$ if $\sigma = \tau$ and $0$ otherwise.
Extend, by linearity, $\inner{\cdot}{\cdot}$ to a scalar product on $\chain{k}$.
Note that we can identify the scalar product $\inner{\cdot}{\cdot}$ with the duality product between chains and cochains in \cref{chains} via the isomorphism $\phi_k : \chain{k} \to \cochain{k}$ in \cref{def_delta}, i.e. $\inner{\sigma}{\tau}=\inner{\phi_k(\sigma)}{\tau}=\phi_k(\sigma)(\tau)$.

Let us consider the chain complex $\ob$ with basis $\basis$.
We define a relation $\prec$ on $\basis$ as follows.
Given distinct basis elements $\sigma \in \basis_k$ and $\tau \in \basis_{k+1}$,
\begin{equation}
\sigma \prec \tau \iff \inner{\sigma}{\bd{k+1} \tau}\neq 0.
\end{equation}
If $\sigma \prec \tau$, then we say that $\sigma$ and $\tau$ are \emph{incident}.

We introduce the \emph{boundary set} and \emph{coboundary set} of $\sigma \in \basis_k$ as
\begin{equation}
\bbd{\basis}{\sigma} \coloneqq \{\,\rho \in \basis_{k-1} \mid \rho \prec \sigma\},
\end{equation}
and
\begin{equation}
\bcobd{\basis}{\sigma} \coloneqq \{\,\rho \in \basis_{k+1} \mid \sigma \prec \rho \},
\label{cobd}
\end{equation}
respectively.

Let $\sigma \in \basis_k$.
If the cardinality of $\bcobd{\basis}{\sigma}$ is one, then we say that $\sigma$ is \emph{free}.
In this case, there exists a unique basis element $\tau$ such that $\sigma \prec \tau$, and we also say that the pair $(\sigma, \tau)$ is free.
If the cardinality of $\bcobd{\basis}{\sigma}$ is greater than one, then we say that $\sigma$ in \emph{internal}.
In this case, if $\tau \in \bcobd{\basis}{\sigma}$, then we also say that the pair $(\sigma,\tau)$ is internal.

\begin{definition}[Matching, acyclic matching]
\label{matching}
A \emph{matching} $\mc M$ on $\basis$ is a family of pairs $\{(\sigma,\tau)\}$ with $\sigma, \tau \in \basis$ such that:
\begin{enumerate}
\item $(\sigma,\tau) \in \mc M$ implies $\sigma \prec \tau$.
\item each $\sigma \in \basis$ is the first component of at most one pair $(\sigma,\tau)$ in $\mc M$.
\end{enumerate}

A matching $\mc M$ is called \emph{acyclic} if there does not exist a cycle
\begin{equation}
\tau_1 \succ \sigma_1 \prec \tau_2 \succ \cdots \prec \tau_h \succ \sigma_h \prec \tau_1,
\label{cycle}
\end{equation}
with $h\geq 2$, $(\sigma_i,\tau_i)\in \mc M$ for all $i \in \{1,\dots, h\}$ and all $\tau_i \in \basis$ being distinct.

A matching $\mc M_k$ of $k$-chains on $\basis$ is a matching such that if $(\sigma, \tau) \in \mc M_k$ then $\sigma \in \basis_k$.
\end{definition}

The following result is a reformulation of Theorem 11.2 in \cite{Kozlov2008} by Kozlov.
It describes the crucial combinatorial property that characterizes acyclic matchings.
Its proof can be obtained by a suitable adaption of the mentioned Theorem 11.2, see pages 181-182 of \cite{Kozlov2008}.
\begin{theorem}
\label{acyclic}
A matching  $\mc M$ on $\mc B$ is acyclic if and only if there exists a total order of pairs of $\mc M$ as $\{(\sigma_1,\tau_1), \dots, (\sigma_n, \tau_n)\}$ such that, for every $i\in\{1,\dots,n-1\}$, $\sigma_i$ is not incident to any $\tau_{i+1},\dots,\tau_n$.
\end{theorem}
In what follows, we will write an acyclic matching $\mc M$ as a sequence $\{(\sigma_1,\tau_1), \dots, (\sigma_n,\tau_n)\}$, where it is understood that the total order is chosen according to \cref{acyclic}.

Let $\mc M$ be a matching on $\basis$.
We say that a basis element in $\basis$ is \emph{matched} (with respect to $\mc M$) if it is contained in some pair in $\mc M$ (both as first or second component of the pair), otherwise, we say that it is \emph{unmatched} or \emph{critical} (with respect to $\mc M$).

We denote by $\mc U_k$ the set of all $\tau \in \basis_k$ such that $\tau$ is matched with some $\sigma \in \basis_{k-1}$.
Similarly, we denote by $\mc D_k$ the set of all $\sigma \in \basis_k$ such that $\sigma$ is matched with some $\tau \in \basis_{k+1}$.
Given $\mc D_k$ and $\mc U_k$, there is a corresponding set of \emph{critical} $k$-chains
\begin{equation}
\mc C_k \coloneqq \basis_k \setminus (\mc D_k \cup \mc U_k).
\end{equation}
Finally, we set $\mc U \coloneqq \bigcup_k \mc U_k$, $\mc D \coloneqq \bigcup_k \mc D_k$ and $\mc C \coloneqq \bigcup_k \mc C_k$.
It easy to see that the sets $\mc U,\mc D, \mc C$ provide a partition of $\basis$
\begin{equation}
\basis = \mc U \cup \mc D \cup \mc C.
\label{partition}
\end{equation}

Our definition of a matching is related to the presentation of the combinatorial Morse theory of Forman \cite{Forman1998MorseTF} and in particular the more recent formulation given by Kozlov \cite{Kozlov2008}.
In earlier presentations, elements in $\mc U$ and $\mc D$ are not explicitly introduced since what is important is only the bijective pairing between their elements.
Instead, in our setting they will play a fundamental role since we use discrete Morse theory from a purely combinatorial point of view and elements in $\mc U$ and $\mc D$ will be used to select suitable submatrices.
The set of critical elements is present also in classical discrete Morse theory.
The set of critical elements will play a fundamental role in \cref{recursive_algo}, where critical elements become the new input for subsequent iterations.
In \cref{recursive_algo}, it is also  essential to be able to express new basis elements as a function of the previous ones.
Hence, we need to keep track of basis structure at each iteration.

An important difference with classical discrete Morse theory is that in Forman \cite{Forman1998MorseTF} a new chain complex, the so-called Morse complex, is constructed.
Our version of discrete Morse theory operates directly on basis elements by performing elementary operations on it.
This is needed to describe how incidence matrices transform in the new bases due to our combinatorial operations.

\subsection{Basis transformations associated with an acyclic matching}
\label{change_basis}
\begin{figure}[h]
    \centering
    \includegraphics[scale=0.8]{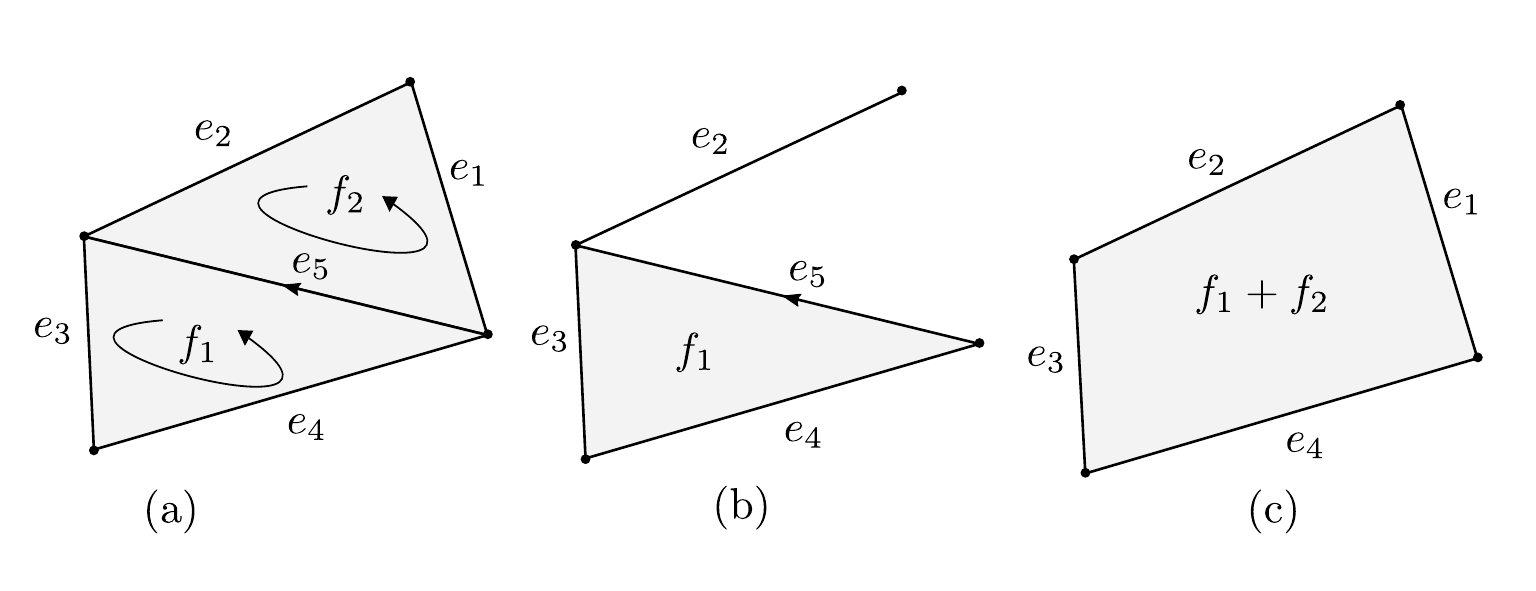}
    \caption{(a) The simplicial complex $\mc K$. (b) Elementary collapse of the free pair $(e_1,f_2)$. (c) Internal collapse of the pair $(e_5,f_2)$; note that the geometric realization of the resulting cell complex is not more simplicial.}
    \label{collapse}
\end{figure}

We begin by considering simple examples to develop geometric intuition behind the general definitions.
When simplifying a cell complex, the effect of a collapse is that of changing the structure of the basis $\basis$, by performing elementary operations on it.
There are three elementary operations to obtain a new basis from a previous one.
If $\basis_k = \{\dots,\xi_i,\xi_j,\dots\}$ is a basis of $\chain{k}$, then a new basis may be obtained by
\begin{enumerate}
\item Exchanging elements $\xi_i$ and $\xi_j$.
\item Multiplying $\xi_i$ by $-1$.
\item Replacing $\xi_j$ by $\xi_j+ q\xi_i$ with $q \in \R$.
\end{enumerate}
Let us consider a prototype example of a simplicial complex $\mc K$ in \cref{collapse}.
The set of all 1-chains is generated by the canonical basis $\Cbasis_1 = \{e_1, e_2, e_3, e_4, e_5 \}$ and the set of all 2-chains by the canonical basis $\Cbasis_2 = \{f_1,f_2\}$.
In \cref{collapse} edges $e_1, e_2, e_3, e_4$ are free since each of them is incident to exactly one face in $\mc K$.
Instead, edge $e_5$ is internal since is the common edge of $f_1$ and $f_2$.
Depending on whether a collapse is elementary or internal we have corresponding elementary operations on $\basis$.

\begin{example}[Elementary collapse]
Let us consider the elementary collapse of the free pair $(e_1,f_2)$.
The obtained cell complex in \cref{collapse}(b) is generated by the set of critical basis elements $\mc C$.
In fact, we have $\mc C_1 = \{e_2,e_3,e_4,e_5\}$ and $\mc C_2=\{f_1\}$.
We get a new basis $\basis'$ of $\mc K$ as $\basis'_1 = \mc D_1 \cup \mc C_1$ and $\basis'_2 = \mc U_2 \cup \mc C_2$, where $\mc D_1 = \{e_1\}$ and $\mc U_2 = \{f_2\}$.
We see that, in the case of an elementary collapse, we get a new basis $\basis'$ of $\mc K$ by performing elementary operations of type 1.
\end{example}

\begin{example}[Interior collapse]
Let us now consider the collapse of the pair $(e_5,f_2)$.
Contrary to the previous case, edge $e_5$ is not free, so we cannot consider an elementary collapse of the pair $(e_5,f_2)$.
However, we can collapse $(e_5,f_2)$ as an internal collapse.
The obtained cell complex in \cref{collapse}(c) is not generated by critical basis elements in $\mc C$ as in the previous case.
In fact, we have $\mc C_1 = \{e_1,e_2,e_3,e_4\}$ and $\mc C_2=\{f_1\}$.
Instead, it is generated by a new set of basis elements obtained from $\mc C$ by adding a linear combination of other basis elements.
We consider the linear transformation $\mc C_2 \ni f_1 \mapsto f_1+f_2$.
The new set of critical basis elements is $\mc C_1 = \{e_1,e_2,e_3,e_4\}$ and $\mc C_2=\{f_1+f_2\}$. 
We get a different basis $\basis'$ of the cell complex $\mc K$ as $\basis'_1 = \mc D_1 \cup \mc C_1$ and $\basis'_2 = \mc U_2 \cup \mc C_2$, where $\mc D_1 = \{e_5\}$ and $\mc U_2 = \{f_2\}$.
We see that, in the case of an internal collapse, we get a new basis $\basis'$ of $\mc K$ by performing elementary operations of type 1, 2 and 3.
\end{example}

Let $\mc M_k = \{(\sigma_1,\tau_1), \dots, (\sigma_n, \tau_n)\}$ be an acyclic matching of $k$-chains on $\basis$.
We now give a recursive definition of the \emph{change of basis associated with $\mc M_k$} and we denote it by $\basis \cdot \mc M_k$.
To start with, given the matched pair $(\sigma,\tau)$, we define the basis $\basis \cdot (\sigma,\tau)$ as follows.
$\basis \cdot (\sigma,\tau)$ is obtained from $\basis$ by performing two actions.
First, we consider the partition of $\basis=\mc U \cup \mc D \cup \mc C$, in particular, we have
\begin{equation}
\basis_k = \mc D_k \cup \mc C_k,
\end{equation}
\begin{equation}
\basis_{k+1} = \mc U_{k+1} \cup \mc C_{k+1}.
\end{equation}
Second, the pair $(\sigma, \tau)$ acts on the set of critical elements $\mc C$ as follows
\begin{equation}
\mc C_{k} \ni \sigma' \mapsto \sigma',
\label{right_action1}
\end{equation}
\begin{equation}
\mc C_{k+1} \ni \tau' \mapsto \tau' - \frac{\inner{\sigma}{\bd{k+1}\tau'}}{\inner{\sigma}{\bd{k+1}\tau}} \tau.
\label{right_action2}
\end{equation}
We see that $\basis \cdot (\sigma,\tau)$ is again a basis of $C^*$.
Indeed, it is obtained from $\basis$ by adding linear combinations of other basis elements.
We define $\basis \cdot \mc M_k$ recursively by the rule
\begin{align}
\basis \cdot \{(\sigma_1,\tau_1),\dots,(\sigma_n,\tau_n)\}\coloneqq
\big(\basis \cdot \{(\sigma_1,\tau_1),\dots,(\sigma_{n-1},\tau_{n-1})\}\big)\cdot (\sigma_n,\tau_n),
\label{recursive_def}
\end{align}
with $n\geq2$.
We see that \cref{recursive_def} is well-defined.
Indeed, if $\basis \cdot \{(\sigma_1,\tau_1),\dots,(\sigma_{n-1},\tau_{n-1})\}$ is basis, then $\basis \cdot \{(\sigma_1,\tau_1),\dots,(\sigma_{n},\tau_{n})\}$ is obtained by adding linear combinations of basis elements, hence it is a basis.
Moreover, $\mc M_k$ is an acyclic matching on $\basis \cdot \{(\sigma_1,\tau_1),\dots,(\sigma_{n-1},\tau_{n-1})\}$, since transformations \cref{right_action1}, \cref{right_action2} leaves matched elements in $\basis$ invariant.

Note that, when $(\sigma,\tau)$ is a free pair, the transformation \cref{right_action2} leaves elements in $\mc C_{k+1}$ invariant.
Hence, no algebraic operations are required and the new basis is simply a permutation of the previous one.
Instead, when $(\sigma,\tau)$ is internal, we have to apply at least one transformation \cref{right_action2}.

\section{Acyclic matchings and Gaussian elimination}
\label{gaussian}
We shall now head towards an algorithm to reduce the matrix $\mat C$ to a row echelon form by means of elementary operations on the basis $\basis$.
The whole procedure boils down to a standard train of thought used in basic linear algebra.
When the coboundary operator $\cobd{1}$ is given as a finite matrix $\mat C$, the bases and orders are already determined.
However, we can get any other bases by applying elementary row and column operations on the matrix $\mat C$.
Our algorithm will produce a change of bases in such a way that the matrix $\mat C$ in the new bases has an invertible upper triangular submatrix (i.e. a matrix with non-zeros only in its upper triangle and main diagonal) and thus it can be transformed in row echelon form.
Hence, we can fast solve the system by processing the unknown variables in reverse order, a standard process known as back substitution.

\subsection{Acyclic matchings and Gaussian elimination}
\label{reduction}

The crucial observation is the following general novel result.

\begin{lemma}
\label{triangular}
Denote by $\mat D_k$ one among the matrices $\mat G$, $\mat C$ or $\mat D$ for $k$ equal to $0$,$1$ or $2$, respectively.
Let $\mc M_k$ be an acyclic matching of $k$-chains on $\basis$.
Then, $\block{\mat D_k}{\mc U_{k+1} \times \mc D_k}$, the submatrix of $\mat D_k$ induced by $\mc U_{k+1} \times \mc D_k$, is upper triangular and invertible.
\begin{proof}
Since the matching $\mc M_k=\{(\sigma_1,\tau_1),\dots,(\sigma_n,\tau_n)\}$ is acyclic, for every $i \in \{1,\dots,n-1\}$, the basis element $\sigma_i$ is not incident to any basis element $\tau_{i+1},\dots,\tau_n$.
Thus, all the non-zero entries in each column are above the diagonal as the rows and columns are arranged in the total order induced by the matching.
\end{proof}
\end{lemma}

By \cref{triangular}, for every acyclic matching there is a corresponding upper triangular submatrix $\block{\mat D_k}{\mc U_{k+1} \times \mc D_k}$ of matched rows and columns of $\mat D_k$.
However, \cref{triangular} is decisive only when the number of matched pairs in $\mc M_k$ is equal to the rank of $\mat D_k$.
In fact, in this case, $\block{\mat D_k}{\mc U_{k+1} \times \mc D_k}$ is an invertible submatrix of $\mat D_k$ of order equal to the rank $\mat D_k$, hence, we can write a solution of \cref{main} after setting some free variables to zero.
For this reason, we introduce the new concept of complete acyclic matching.
\begin{definition}[Complete acyclic matching]
We say that a matching $\mc M_k$ of $k$-chains on $\basis$ is \emph{complete} if the number of matched pairs in $\mc M_k$ is equal to the rank of $\mat D_k$.
\end{definition}
Now, we illustrate how the upper triangular submatrix $\block{\mat D_k}{\mc U_{k+1} \times \mc D_k}$ induced by a complete acyclic matching allows us to operationally obtain a solution of \cref{main} using back substitution.

Let us consider a complete acyclic matching $\mc M_1$ of 1-chains.
By applying \cref{triangular}, the action of $\mc M_1$ on \cref{main} is equivalent to Gaussian elimination.
It is thus sufficient to invert submatrix $\block{\mat C}{\mc U_2 \times \mc D_1}$ after setting free variables in $\mc C_1$ to zero.
We write a discrete potential $\bm h$ solution of \cref{main} as
\begin{align}
&\block{\bm h}{\mc D_1} = \block{\mat C}{\mc U_2 \times \mc D_1}^{-1} \block{\bm i}{\mc U_2},\\
&\block{\bm h}{\mc C_1} = \bm 0.
\end{align}
Since $\block{\mat C}{\mc U_2 \times \mc D_1}$ is upper triangular, we can evaluate  $\block{\mat C}{\mc U_2 \times \mc D_1}^{-1} \block{\bm i}{\mc U_2}$ by back substitution in linear time \cite{Strang1993IntroductionTL}.

The aim of the back substitution is to determine the coefficient values of $\bm h$.
Let $\mc M_1 = \{(e_1, f_1),\dots,(e_n,f_n)\}$.
The process of back substitution is so-called because one determines the coefficient values backwards, by first computing $\bm h_{e_n}$, then substituting back into the previous equation to solve for $\bm h_{e_{n-1}}$ and repeating through $\bm h_{e_{n-1}}$.
A naive combination of these coefficients in one step for each $\bm h_{e_i}$ is very time consuming since large intermediate expressions are generated.
This is avoided by combining coefficients pairwise as in the following standard back substitution algorithm.
\begin{algorithm}[H]
\caption{Back substitution process}
\begin{algorithmic}[1]
\Procedure{Backsubstitution}{$\mc M_1, \basis$}
	\For{$i$ ranging from $n$ down to $1$}
		\State $\bm h_{e_i} = (\mat C_{f_i,e_i})^{-1} (\bm i_{f_i} - \sum_{e \in \bbd{\basis}{f_i}} \mat C_{f_i,e} \bm h_e)$
	\EndFor
	\State \textbf{return} $\bm h$;
\EndProcedure
\end{algorithmic}
\label{back_algo}
\end{algorithm}

In the above discussion, although we focused on matrix $\mat C$, we have actually detailed the proof of the following theorem.

\begin{theorem}
\label{main_lemma}
Given a complete acyclic matching $\mc M_k$ of $k$-chains on $\basis$, for $k \in \{0,1,2\}$ we can find in linear time a discrete potential $\bm v \in \cochain{k}$ of $\bm w \in \cochain{k+1}$ with $\mat D_{k+1}\bm w = \bm 0$, namely a solution of $\mat D_k \bm v = \bm w$.
\end{theorem}

\subsection{On the problem of constructing complete acyclic matchings: the case of tree-cotree techniques}
\label{acyclic_matching_sec}
We now turn to the theoretical issue of constructing complete acyclic matchings $\mc M_k$ of $k$-chains on the canonical basis $\Cbasis$.

Case $k=0$.
We need to construct a complete acyclic matching $\mc M_0$ of 0-chains on $\Cbasis$.
By definition, the number of matched elements in $\mc M_0$ has to be equal to the rank of $\mat G$, which is $\vc -1$.
Let us consider a \emph{spanning tree} $T$ on $\mc K$.
It can be constructed in worst-case linear time using standard graph algorithms, for instance \emph{breadth-first search} (BFS) algorithm \cite{Cormen2009IntroductionTA}.
There is a standard reasoning to define an acyclic matching corresponding to a spanning tree $T$ \cite{Lewiner2003OptimalDM}.
It is obtained by mimicking a spanning tree traversal process.
The construction goes as follows.
Since $T$ is a tree, there exists at least a \emph{leaf} in $T$, namely a vertex $v \in \mc K_0$ with only one incident edge in $T$.
Pick a leaf $v$ and pair it with the unique edge $e$ containing it.
Then, add the pair $(v,e)$ to $\mc M_0$ and remove $v,e$ from $T$.
Since $v$ is a leaf, the obtained graph is again a tree.
By iteratively repeating this process, we define a complete acyclic matching $\mc M_0$.
The matching $\mc M_0$ is well-defined: if a vertex $v$ is matched during the process, it is removed from $T$ and thus cannot appear in any other pair; it is acyclic since every tree is.
The matching $\mc M_0$ is complete, since $\mc K$ is connected, $T$ is spanning, namely all vertices of $\mc K$ are in $T$ and they are eventually added to $\mc M_0$ in the above process, except the last one.

Case $k=2$.
Before considering the case $k=1$, we show the similarity of the construction for $k=2$ with the case $k=0$.
We need to construct a complete acyclic matching $\mc M_2$ of 2-chains on basis $\Cbasis$.
By definition, the number of matched elements in $\mc M_2$ has to equal to the rank of $\mat D$, which is $\cc$.
Note that the cell complex $\mc K$ defines a manifold with boundary.
It implies that if $c,c'$ are distinct volumes and $f$ is a face such that $f \subset \partial c \cap \partial c'$, then $c, c'$ are the only volumes that contain $f$.
Thanks to the manifold condition, it is well-defined the so-called \emph{complete dual graph} of $\mc K$.
The dual graph $\dual G = (\dual V \cup \dual V_\infty, \dual E \cup \dual E_\infty)$ of $\mc K$ is a graph with set of vertices $\dual V$ given by volumes of $\mc K$ and $\{c,c'\}\in \dual E$ if $\dim (c \cap c')=2$.
There is an additional vertex $\dual v_\infty \in \dual V_\infty$ and there are additional edges $\dual E_\infty$ in $\dual G$.
Each edge in $\dual E_\infty$ corresponds to a volume $c$ whose boundary $\partial c$ contains a boundary face $f$; thus, if $f$ is a boundary face and $c$ is the unique volume incident to it, then $\{\dual v_\infty, c\} \in \dual E_\infty$.
Let us consider a spanning tree $\dual{T}$ on $\dual G$. 
By repeating the same construction detailed  for $k=0$ on $\dual T$, we obtain a finite sequence $S$ of pairs of the form $(c,\dual e)$, where $c$ is a vertex in $\dual V \cup \dual V_\infty$ (i.e. a volume of $\mc K$ or $\dual v_\infty$) and $\dual e=\{c,c'\}$ is an edge in $\dual E \cup \dual E_\infty$.
For each pair $(c,\dual e)$, we choose a face $f$ of $\mc K$ such that $f \subset \partial c \cap \partial c'$ if $\dual v_\infty \notin \dual e=\{c,c'\}$, $f \subset \partial c$ if $c'=\dual v_\infty$ and $f \subset \partial c'$ if $c=\dual v_\infty$.
Replace in $S$ each pair $(c,\dual e)$ with either $(f,c)$ if $c \neq \dual v_\infty$ or $(f,c')$ if $c=\dual v_\infty$.
In this way, we obtain a complete acyclic matching $\mc M_2$ of 2-chains of $\mc K$. 
Thus, all vertices of the dual graph are eventually added to $\mc M_2$, except the vertex $\dual v_\infty$, and all volumes in $\mc K$ are matched.

By combining \cref{main_lemma} with the above results, we state the following theorem which solves the discrete potential problem for $k\in\{0,2\}$.

\begin{theorem}
For $k\in\{0,2\}$, we can find in linear time a discrete potential $\bm v \in \cochain{k}$ of $\bm w \in \cochain{k+1}$ with $\mat D_{k+1} \bm w = \bm 0$, namely a solution of $\mat D_k \bm v = \bm w$.
\end{theorem}

Case $k=1$.
With no surprise, it turns out to be the most challenging case.
To begin with, we show the following result.
\begin{theorem}
\label{np_complete}
Let $\mc K$ be a 3-dimensional topologically trivial simplicial complex embedded in $\R^3$ with Lipschitz boundary as in \cref{chains}.
Then, $\mc K$ is collapsible if and only if there exists a complete acyclic matching $\mc M_1$ of 1-chains on the canonical basis $\Cbasis$.
\begin{proof}
Let us consider a sequence of elementary collapses $\mc M$ of $\mc K$ leading from $\mc K$ to a vertex.
We define an acyclic matching $\mc M_1$ of 1-chains (on $\Cbasis$) by selecting, according to the total order of $\mc M$, all pairs made by matched edges and faces in $\mc M$.
Note that $\mc M_1$ is necessarily acyclic since $\mc M$ is.
Next, since $\mc K$ collapses to a vertex, all faces in $\mc K$ are matched in $\mc M$.
In particular, all volumes in $\mc K$ are matched in $\mc M$ and there are $\cc$ corresponding matched faces.
Thus, there are $\fc - \cc$ matched faces in $\mc M_1$.
It follows that $\mc M_1$ is complete since the rank of $\mat C$ is exactly $\fc-\cc$.

To prove the converse result, let us consider a complete acyclic matching of 1-chains $\mc M_1$.

We need the following construction.
First let us consider the complete dual graph $\dual G = (\dual V \cup \dual V_\infty, \dual E \cup \dual E_\infty)$ of $\mc K$ as defined above for the case $k=2$.
Consider the subgraph $\dual G_{\mc M_1}=(\dual V \cup \dual V_\infty, \dual E_{\mc M_1})$ of $\dual G$ which has the same vertices of $\dual G$ and $\dual E_{\mc M_1}$ includes an edge $\{c,c'\} \in \dual E$ of $\dual G$ if the unique face $f$ of $\mc K$ such that $f \subset c \cap c'$ is critical with respect to $\mc M_1$.
We now show that $\dual G_{\mc M_1}$ has $\cc$ edges (i.e, the number of vertices of $\dual G$ minus one) and is connected.
First, since $\mc M_1$ is a complete acyclic matching the number of critical faces, which is equal to the number of edges in $\dual E_{\mc M_1}$, is $\fc - (\fc - \cc )= \cc$.
Second, assume that, for the sake of contradiction, $\dual G_{\mc M_1}$ is disconnected.
Let $\dual V^*$ be the set of nodes in a fixed connected component of $\dual G_{\mc M_1}$.
Let $\dual E_{\dual V^*} \subset \dual E \cup \dual E_\infty$ be the subset of \emph{cut edges}, namely edges of $\dual G$ with one vertex in $\dual V^*$ and one vertex in its complement $(\dual V \cup \dual V_\infty) \setminus \dual V^*$.
Note that since $\dual G$ is connected, $\dual E_{\dual V^*}$ is not empty.
Moreover, by definition of $\dual G_{\mc M_1}$, to each $\dual e \in \dual E_{V^*}$ corresponds a matched face $f$ of $\mc K$ with respect to $\mc M_1$ forming the pair $(e,\dual e)=(e,f) \in \mc M_1$ for a unique edge $e$ of $\mc K$.

Now we construct a cycle of the form $\dual e_1 \succ e_1 \prec \dual e_2 \succ \dots \prec \dual e_h \succ e_h \prec \dual e_1$ with $h\geq2$, $(e_i,\dual e_i) \in \mc M_1$ for all $i \in \{1,\dots,h\}$ and all $\dual e_i \in \dual E_{V^*}$ being distinct.
Start with a cut edge $\dual e_1$ and consider the unique edge $e_1$ of $\mc K$ forming the matched pair $(e_1,\dual e_1) \in \mc M_1$.
Then, $e_1$ is contained in at least one other cut edge $\dual e_2 \in \dual E_{V^*}$, otherwise $\dual e_1$ cannot be a cut edge.
By iteratively repeating this process, we obtain a sequence $S$ of pairs $(e_1,\dual e_1), \dots, (e_i,\dual e_i)$ such that $\dual e_1 \succ e_1 \prec \dual e_2 \succ \dots \prec \dual e_{i-1} \succ e_i \prec \dual e_i$.
Since we have finite graphs, eventually we will run out of those cut edges which do not appear as the second component of a pair in $S$.
Thus, at some step $i$ of the above process, we must get a cut edge that is the second component of a pair in $S$ and we get a cycle of the form \cref{cycle}, which is a contradiction since $\mc M_1$ is acyclic.
We have proved that $\dual G_{\mc M_1}$ is a connected subgraph of $\dual G$ with $\cc$ edges and thus it is a spanning tree of $\dual G$.
We define a complete acyclic matching $\mc M_2$ of 2-chains by using the spanning tree construction detailed above for $k=2$ on $\dual G_{\mc M_1}$.

Now, $\mc M_2$ and $\mc M_1$ define a sequence of elementary collapses leading from $\mc K$ to a single vertex as follows.
Since every volume of $\mc K$ is matched with respect to $\mc M_2$, the set of critical cells of $\mc K$ with respect to $\mc M_2$ form a 2-dimensional subcomplex $\mc K^{\,(1)}$ of $\mc K$.
By applying Theorem 11.13 (a) in \cite{Kozlov2008}, there exists a sequence of elementary collapses leading from $\mc K$ to $\mc K^{\,(1)}$.
Next, observe that, by construction, the set of matched faces in $\mc K$ with respect to $\mc M_2$ is exactly the set of critical faces of $\mc K$ with respect to $\mc M_1$.
Thus, the set of critical cells of $\mc K^{\,(1)}$ with respect to $\mc M_1$ form a 1-dimensional subcomplex $\mc K^{\,(2)}$ of $\mc K^{\,(1)}$. 
By applying Theorem 11.13 (a) in \cite{Kozlov2008} there exists a sequence of elementary collapses leading from $\mc K^{\,(1)}$ to $\mc K^{\,(2)}$.
$\mc K^{\,(2)}$ is a topologically trivial 1-dimensional subcomplex, thus it is a spanning tree on $\mc K$.
Thus $\mc K$ collapses to a single vertex.
\end{proof}
\end{theorem}

Now we focus on the \emph{Spanning Tree Technique} (STT) algorithm introduced in \cite{Dlotko2010CriticalAO} and implicitly used in many works \cite{clenet1998,sinum1,Rodrguez2015FiniteEP,Webb1989ASS}.
The algorithm, given as input a spanning tree $T$ on $\mc K$, computes a discrete vector potential $\bm h$ as follows.
To start with, STT sets the value $\bm h_e = 0$ if $e\in T$, otherwise the value $\bm h_e$ at this stage is unknown.
Next, all faces $f$ of $\mc K$ are loaded into a list $L$.
The main loop of STT works until there are no more faces in $L$.
In each iteration, we randomly search for a face $f$ with two boundary edges $e_1,e_2 \subset \partial f$ such that the values $\bm h_{e_1}, \bm h_{e_2}$ are known. Then, the value of $\bm h_{e}$ of the remaining boundary edge $e \subset \partial f$  is determined by
\begin{equation}
\bm h_{e} = (\mat C_{f,e})^{-1}(\bm i_{f} - \mat C_{f,e_1} \bm h_{e_1} - \mat C_{f,e_2} \bm h_{e_2})
\label{use}
\end{equation}
and face $f$ is removed from the list $L$.
In the case when $L$ is non-empty and there is no available face $f$ satisfying the above property, then STT does not terminate since it stalls in a infinite loop.
As shown in \cite{Dlotko2010CriticalAO}, STT termination depends on the choice of the input spanning tree $T$.

We now show that the STT algorithm boils down to a procedure to construct complete acyclic matchings of 1-chains on the canonical basis $\Cbasis$.
This result completes the picture of the equivalence between tree-cotree techniques and the problem of finding complete acyclic matching of $k$-chains.
The idea behind this observation has it root on the fact that tree-cotree techniques describe the same actions of acyclic matchings although using a different language.

To formally state the next theorem we introduce the following notations.
Given a spanning tree $T$ on $\mc K$, we say that STT \emph{terminates} (with input $T$) if it does not stall in an infinite loop.
We say that STT \emph{uses the pair $(e,f)$} if STT determines the value of edge $e \subset \partial f$ via \cref{use} during its main loop execution.
Finally, we say that STT \emph{uses edge $e$} if STT uses the pair $(e,f)$ for some face $f \in L$.

\begin{theorem}
\label{STT_matching}
Let $\mc K$ be a simplicial complex as in \cref{np_complete}.
There exists a spanning tree on $\mc K$ for which STT terminates if and only if there exists a complete acyclic matching of 1-chains on the canonical basis $\Cbasis$. 
\begin{proof}
Let $T$ be a spanning tree such that STT terminates.
We now define an acyclic matching $\mc M_T$ of 1-chains on $\Cbasis$ such that every \emph{cotree edge} of $T$, namely an edge $e \in \mc K_1$ which is not an edge of $T$, is matched in $\mc M_T$.
$\mc M_T$ is constructed during STT execution as follows.
We initialize $\mc M_1=\emptyset$.
Next, if STT uses a pair $(e,f)$ during its  execution then we add $(e,f)$ to $\mc M_1$.
In this way we get a total order of pairs in $\mc M_1$ where a pair $(e_i,f_i) \in \mc M_1$ comes before than a pair $(e_j,f_j) \in \mc M_1$ in this total order with $i < j$ if STT uses $(e_i,f_i)$ before than $(e_j,f_j)$.
The set $\mc M_T$ of all such pairs is a complete acyclic matching.
First, we see that $\mc M_T$ is a matching since if STT uses the pair $(e,f)$, then the value of $\bm h_e$ set by STT can never be reassigned.
Second, the matching $\mc M_T$ is acyclic since if STT uses the pair $(e_i,f_i)$, then $e_i$ is not contained in the boundary of any matched face $f_j$ in $\mc M_1$ for $1\leq j < i$.
Third, since STT terminates, every cotree edge $e$ of $T$ is matched in $\mc M_1$.
Using Euler's formula, the fact that $\mc K$ is topologically trivial and that $T$ is a spanning tree on $\mc K$ we get $\ec - (\vc - 1) = \fc -\cc$. 
Thus $\mc M_1$ is complete.

Conversely, let $\mc M_1$ be a complete acyclic matching of 1-chains on $\Cbasis$.
We need the following construction, which mimics the one described in \cref{np_complete}.
Consider the subgraph $G_{\mc M_1}=(V, E_{\mc M_1})$ of $\mc K$ with set of vertices given by vertices of $\mc K$ and $E_{\mc M_1}$ includes an edge $e \in \mc K_1$ if $e$ is critical with respect to $\mc M_1$.
We now show that $G_{\mc M_1}$ has $\vc-1$ edges and is connected.
First, since $\mc M_1$ is a complete acyclic matching the number of critical edges, which is equal to the number of edges in $E_{\mc M_1}$, is $\ec - (\fc - \cc )= \vc -1$ where we have used Euler's formula $\vc - \ec +\fc -\cc = 1$ and the fact that $\mc K$ is topologically trivial.
Second, assume that, for the sake of contradiction, $G_{\mc M_1}$ is disconnected.
Let $V^*$ be the set of nodes in a fixed connected component of $G_{\mc M_1}$.
Let $E_{V^*} \subset \mc K_1$ be the subset of \emph{cut edges}, namely edges of $G$ with one vertex in $V^*$ and one vertex in its complement $V \setminus V^*$.
Note that since $\mc K$ is connected, $E_{V^*}$ is not empty.
Moreover, by definition of $G_{\mc M_1}$, to each $e \in E_{V^*}$ corresponds a matched edge $e$ of $\mc K$ with respect to $\mc M_1$ forming the pair $(e,f) \in \mc M_1$ for a unique face $f$ of $\mc K$.

Now we construct a cycle of the form $f_1 \succ e_1 \prec f_2 \succ \dots \prec f_h \succ e_h \prec f_1$ with $h\geq2$, $(e_i, f_i) \in \mc M_1$ for all $i \in \{1,\dots,h\}$ and all $e_i \in E_{V^*}$ being distinct.
Start with a cut edge $e_1$ and consider the unique face $f_1$ of $\mc K$ forming the matched pair $(e_1,f_1) \in \mc M_1$.
Then, $f_1$ contains at least one other cut edge $e_2 \in E_{V^*}$, otherwise $e_1$ cannot be a cut edge.
By iteratively repeating this process, we obtain a sequence $S$ of pairs $(e_1,f_1), \dots, (e_i,f_i)$ such that $f_1 \succ e_1 \prec f_2 \succ \dots \prec f_{i-1} \succ e_i \prec f_i$.
Since we have finite graphs, eventually we will run out of those cut edges which do not appear as the first component of a pair in $S$.
Thus, at some step $i$ of the above process, we must get a cut edge that is the first component of a pair in $S$ and we get a cycle of the form \cref{cycle}, which is a contradiction since $\mc M_1$ is acyclic.
We have proved that $G_{\mc M_1}$ is a connected subgraph of $\mc K$ with $\vc-1$ edges and thus it is a spanning tree on $\mc K$.

Now we prove that STT with input $G_{\mc M_1}$ terminates.
Note that, by definition of $G_{\mc M_1}$, if $(e,f) \in \mc M_1$ then $e$ is a \emph{cotree edge} of $G_{\mc M_1}$, namely $e$ is not an edge of $G_{\mc M_1}$.
Denote by $C$ the set of all cotree edges of $G_{\mc M_1}$.
Moreover, since $\mc M_1$ is complete, reasoning as above using Euler's formula, we get that every cotree edge of $G_{\mc M_1}$ is matched with respect to $\mc M_1$.
Since $\mc M_1=\{(e_1,f_1),\dots,(e_n,f_n)\}$ is acyclic, we can order pairs in $\mc M_1$ in such a way that, for every $i \in \{1,\dots, n-1\}$, $f_i$ is not incident to any $e_{i+1}, \dots,e_{n}$.
Let $E_{\text{STT}}$ be the set of cotree edges used by STT during its execution.
Note that $E_{\text{STT}}$ is not empty since STT can use at least the pair $(e_1,f_1)$, thanks to the total order chosen on $\mc M_1$.
Suppose that STT does not terminate.
This means that, during its execution, it does not use any cotree edge and $C \setminus E_{\text{STT}} \neq \emptyset$.
We will show that this is impossible.
Let $i$ be the minimum integer such that the cotree edge $e_i$ belongs to $C \setminus E_{\text{STT}}$.
Since pairs in $\mc M_1$ are ordered as described above, $i >1$ and $f_i$ can be only incident to cotree edges $e_j$ with $j \leq i$.
Moreover, by definition of $i$, there is no $j$ with $j < i$ such that the cotree edge $e_j$ belongs to $C \setminus E_{\text{STT}}$.
Hence, STT should have at least used the pair $(e_i,f_i)$, i.e. $e_i$ belongs to $E_{\text{STT}}$ as well.
This gives the desired contradiction and completes the proof.
\end{proof}
\end{theorem}

The proof of \cref{STT_matching} shows that, if the STT algorithm terminates for a given spanning tree $T$, starting from $T$ we can construct a complete acylclic matching $\mc M_T$ of 1-chains on $\Cbasis$.
But also the other way around, namely, if we have a complete acyclic matching $\mc M_1$ of 1-chains on $\Cbasis$, starting from $\mc M_1$ we can construct a spanning tree $G_{\mc M_1}$ for which STT terminates.

By combining \cref{np_complete} and \cref{STT_matching} we now state the following result which gives a topological characterization of termination problems of tree-cotree techniques.

\begin{theorem}
\label{characterization}
Let $\mc K$ be a simplicial complex as in \cref{np_complete}.
Then, there exists a spanning tree on $\mc K$ for which STT terminates if and only if $\mc K$ is collapsible.
\end{theorem}

There are known examples of triangulations of 3-balls which are not collapsible \cite{Lutz2013KnotsIC}.
Hence, using \cref{characterization}, there are triangulations $\mc K$ of 3-balls such that, for every possible spanning tree of $\mc K$ given as input, STT does not terminate.

We remark that it is NP-complete to decide whether a given 3-dimensional simplicial complex (embedded or not) is collapsible \cite{Tancer2016RecognitionOC}.
However, to the authors knowledge, the related question for the case of 3-dimensional simplicial complexes embedded in $\R^3$ and with Lipschitz boundary is still an open problem. 
Yet, numerical evidence shows that this problem is very difficult in general although good heuristics exists \cite{Benedetti2014RandomDM}.

These results are the cause of the well-known termination problems of tree-cotree techniques \cite{Dlotko2010CriticalAO}; see also the discussion in \cref{numerical}.

\subsection{The case where we cannot find a complete acyclic matching}
\label{residual_sec}
As shown by \cref{np_complete} and the discussion at the end of \cref{acyclic_matching_sec}, for cell complex $\mc K$ given as input, we cannot find in general a complete acyclic matching $\mc M_1$ of 1-chains on the canonical basis $\Cbasis$.
Consequently, \cref{main_lemma} cannot be applied in general for $\basis = \Cbasis$.

Let us now assume that acyclic matching $\mc M_1$ is not complete.
We shall now show that, after solving another linear system, we can still get a discrete vector potential solution $\bm h$ solution of \cref{main} by exploiting back substitution.

As shown in \cref{change_basis}, an acyclic matching acts on a basis $\basis$ by performing elementary operations on it.
Let us consider the new basis $\basis' = \basis \cdot \mc M_1=\mc U \cup \mc D \cup \mc C$.
There is a corresponding block partition of linear system \cref{main} as
\begin{equation}
\begin{pmatrix}
\block{\mat C}{\mc U_2 \times \mc D_1} & \block{\mat C}{\mc U_2 \times \mc C_1}\\
\block{\mat C}{\mc C_2 \times \mc D_1} &\block{\mat C}{\mc C_2 \times \mc C_1} 
\end{pmatrix}
\begin{pmatrix}
\block{\bm h}{\mc D_1}\\
\block{\bm h}{\mc C_1}
\end{pmatrix}
=
\begin{pmatrix}
\block{\bm i}{\mc U_2}\\
\block{\bm i}{\mc C_2}
\end{pmatrix}.
\end{equation}
The crucial fact turns out to be that $\block{\mat C}{\mc C_2 \times \mc D_1}$ is a zero matrix.
Thus, we can determine a discrete vector potential $\bm h=(\block{\bm h}{\mc D_1},\block{\bm h}{\mc C_1})^T$ solution of \cref{main} by solving, in order, 
\begin{equation}
\block{\mat C}{\mc C_2 \times \mc C_1} \block{\bm h}{\mc C_1} = \block{\bm i}{\mc C_2},
\label{residual}
\end{equation}
\begin{equation}
\block{\mat C}{\mc U_2 \times \mc D_1} \block{\bm h}{\mc D_1} = \block{\bm i}{\mc U_2} - \block{\mat C}{\mc U_2 \times \mc C_1} \block{\bm h}{\mc C_1},
\label{residual2}
\end{equation}
where \cref{residual2} is solved by exploiting back substitution as in \cref{reduction}.

Let us prove that $\block{\mat C}{\mc C_2 \times \mc D_1}$ is a zero matrix.
Let $\mc M_1=\{(\sigma_1, \tau_1) \dots, (\sigma_n,\tau_n)\}$.
Let us consider bases $\{\dots,\xi^i=\phi_k(\xi_i),\dots\}$ and $\{,\dots,\pi_i,\dots,\}$ of $\cochain{1}$ and $\chain{2}$, respectively.
We see that the $(i,j)$-entry of $\mat C$ is $\mat C_{i,j} = \inner{\cobd{1}(\xi^j)}{\pi_i}=\inner{\xi^j}{\bd{2}\pi_i}=\inner{\xi_j}{\bd{2}\pi_i}$, where we have used \cref{def_cobd} and the isomorphism $\phi_k : \chain{k} \to \cochain{k}$ in \cref{def_delta}.
Thus, to prove that $\block{\mat C}{\mc C_2 \times \mc D_1}$ is a zero matrix we have to show that every $\sigma_i \in \mc D_1$ with $i \in \{1,\dots,n\}$ is not incident to any $\rho \in \mc C_2$.

We proceed by induction on $n$.

Suppose that $n=1$.
Necessarily, $\tau_1 \in \mc U_2$ is the only basis element in $\basis\cdot \mc M_1=\mc U \cup \mc D \cup \mc C$ incident on $\sigma_1$.
In fact, if $\tau \in \basis$ is incident on $\sigma_1$, then its image under the transformation \cref{right_action2} is not incident on $\sigma_1$.
It follows that $\sigma_1$ is not incident to any $\rho \in \mc C_2$.

We now assume the statement true for $n-1$ and we prove it for $n$.
Let us consider the basis $\basis \cdot \mc M_1 = \mc U \cup \mc D \cup \mc C$.
Proceeding as above, we see that $\sigma_n$ is not incident to any $\rho \in \mc C_2$.
To conclude, it is sufficient to show that also each $\sigma_i$ with $i\in\{1,\dots,n-1\}$ is not incident to any $\rho \in \mc C_2$.
Basis elements in $\basis \cdot \mc M_1$ are obtained from that of $\basis \cdot \{(\sigma_1,\tau_1),\dots,(\sigma_{n-1},\tau_{n-1})\}$ by applying transformation \cref{right_action2} with $(\sigma,\tau)=(\sigma_n,\tau_n)$.
Using the induction hypothesis we see that the first term in \cref{right_action2} (i.e., $\tau' \in \chain{2}$) is not incident on any $\sigma_i$ with $i \in\{1,\dots,n-1\}$.
Moreover, each $\sigma_i$ with $i \in \{1,\dots,n-1\}$ is not incident to $\tau_n$ as pairs in $\mc M_1$ are ordered as in \cref{acyclic}.
This completes the proof.

We now state the following theorem which combines the results of this section with those of \cref{reduction}.
It can be thought as a specific algebraic version of Forman's discrete Morse complex construction \cite{Forman1998MorseTF}.
We think that our presentation and terminology shed light on the linear algebra behind the more abstract discrete Morse theory constructions.
\begin{theorem}
Denote by $\mat D_k$ one among the matrices $\mat G$, $\mat C$ or $\mat D$ for $k$ equal to $0$,$1$ or $2$, respectively.
Let $\mc M_k$ be an acyclic matching of $k$-chains on $\basis$.
Then, there is a corresponding block partition of $\mat D_k$ as 
\begin{equation}
\mat D_k=
\begin{pmatrix}
\block{\mat D_k}{\mc U_{k+1} \times \mc D_k} & \block{\mat D_k}{\mc U_{k+1} \times \mc C_k}\\
0 &\block{\mat D_k}{\mc C_{k+1} \times \mc C_k} 
\end{pmatrix},
\end{equation}
where $\block{\mat D_k}{\mc U_{k+1} \times \mc D_k}$ is upper triangular and invertible.
\end{theorem}

\section{Algorithm description}
\label{algorithm}
The goal of this section is a recursive algorithm that reduces matrix $\mat C$ into a row echelon form by means of elementary operations on the basis $\basis = \bigcup_k \basis_k$ and we present in \cref{recursive}.
We first describe in \cref{greedy} our novel greedy procedure to construct acyclic matchings. 
\subsection{Greedy approach to construct acyclic matchings}
\label{greedy}
To minimize computational effort of change of basis in \cref{residual}, we have to carefully choose how to construct the acyclic matching $\mc M_1$.
It is visible from \cref{change_basis} that constructing acyclic matchings by internal collapses can get more complicated than by elemetary collapses since changes of basis of type 2 and 3 are involved.
In this case, we have to express the new basis elements as a linear combination of the previous ones and hence they cannot be removed from the data structure after each collapse.

The above discussion motivates the concept of \emph{degree} of a basis element $\sigma \in \basis$.
We define $\degree{\sigma}$ as the cardinality of the coboundary of $\sigma$, i.e. the set $\bcobd{\basis}{\sigma}$ defined in \cref{cobd}.
Note that if $\degree{\sigma}$ is 1 then $\sigma$ is free.
In our approach, also the case $\degree{\sigma} =2$ will play a fundamental role.
We define $\sigma \in \basis$ to be \emph{flat} if $\degree{\sigma}=2$.

To reduce the amount of computation during Gaussian elimination, it is wise to first search for basis elements with lowest degree.
Indeed, if $\degree{\sigma}$ is 1, then $\sigma$ is free and the new basis is a selection of the previous one and no algebraic operations are needed.
In this case, it should be possible to avoid all the matrix algebra computation and efficiently organize the basis so that matrices are in triangular form as in \cref{triangular}.

Keeping in mind the above heuristics, we search for basis elements having smallest degree.
Specifically, we do not strictly choose basis elements with minimum degree but instead we proceed in a sequential manner.
First, we search for all free basis elements until exhaustion.
Next, we search for all flat basis elements until exhaustion.
We have pursued this method because been motivated by its practical implementation and performance on test problems, rather than by following the best theoretical greedy approach; see the discussion in \cref{numerical}.

We construct acyclic matchings using a standard elementary collapse greedy procedure, where we search for collapsing sequences of free basis elements in a monotone-like fashion \cite{Benedetti2014RandomDM}.
That is, we proceed in sequential order with respect to the dimension of the basis elements by first collapsing 2-chains and then 1-chains.
\begin{algorithm}[H]
\caption{Random strategy to construct acyclic matchings}
\begin{algorithmic}[1]
\Procedure{ConstructAcyclicMatchings}{$\basis$}
	\For{$k$ ranging from $2$ down to $1$}
	\While{there exists a free pair $(\sigma,\tau)$ in $\basis_k$}
		\State elementary collapse $(\sigma,\tau)$
		\State insert $(\sigma,\tau)$ to $\mc M_k$
	\EndWhile
	\While{there exists a flat pair $(\sigma,\tau)$ in $\basis_k$}
		\State internal collapse of $(\sigma,\tau)$
	    \State insert $(\sigma,\tau)$ to $\mc M_k$
	\EndWhile
	\EndFor
	\State \textbf{return} $\mc M_2 \cup  \mc M_1$
\EndProcedure
\end{algorithmic}
\label{acyclic_algo}
\end{algorithm}

It is clear that \Cref{acyclic_algo} always terminates.
Moreover, the obtained matching is acyclic since at each iteration collapses are performed.
Note that \Cref{acyclic_algo} requires no backtracking since new free or flat pairs can only appear after each new collapse.
Thus, the worst-case complexity is linear using a suitable algorithm implementation that employs a list data structure.

As proved in \cref{acyclic_matching_sec}, we can always find a complete acyclic matching of 2-chains by using a standard spanning tree construction.
The proof of this result implies that the order in which we collapse 2-chains in \cref{acyclic_algo} is not important.
Therefore, there always exists a complete acyclic matching of 2-chains and we can get one by performing collapses in a random fashion.

\subsection{Recursive algorithm}
\label{recursive}
We now present a recursive construction of acyclic matchings.

Normally one aims at finding a complete acyclic matching, namely a matching that reaches the needed number of matched pairs so that we can apply \cref{main_lemma}.
However, as proved in \cref{np_complete}, determining a complete acyclic matching $\mc M_1$ of 1-chains is a hard algorithmic problem and we do not tackle it.
In \cref{acyclic_algo} we employed a greedy strategy that randomly selects an acyclic matching but it can easily happen that the obtained acyclic matching $\mc M_1$ is not complete.

Let us consider the case where $\mc M_1$ is not complete.
After applying basis transformations associated with $\mc M_1$ to the current basis $\basis$, we get a new basis $\basis' = \basis \cdot \mc M_1=\mc U \cup \mc D \cup \mc C$.
In particular, the set of critical basis elements $\mc C$ in $\basis'$ forms a linearly independent subset of $\basis'$.
As pointed out in \cref{residual_sec}, the focus now shifts to the set of critical basis elements $\mc C$.

The novel idea is that we may call the whole routine recursively, where the output set of critical basis elements obtained from the previous iteration becomes the input basis for the next iteration.
Fundamentally, we think linear system \cref{residual} as an instance of the original linear system \cref{main}, although considering a subset of the previous basis.
This operation is well-defined since we have the partition $\basis' = \mc U \cup \mc D \cup \mc C$ and the set $\mc C$ of critical basis elements form a linearly independent subset of $\basis$.

If at a certain recursion stage we get a complete acyclic matching $\mc M_1$ of 1-chains, then we can find a discrete vector potential $\bm h$ solution of \cref{main} by recursively applying the reasoning described in \cref{residual_sec}.
Otherwise, we recursively apply the same routine on the obtained set of critical basis elements $\mc C$.

\begin{algorithm}[H]
\caption{Construction of a discrete vector potential $\bm h$ solution of \cref{main}}
\begin{algorithmic}[1]
\Procedure{ConstructDiscreteVectorPotential}{$\basis$} \Comment{at the beginning $\basis$ is the canonical basis $\Cbasis$}
	\State construct an acyclic matching $\mc M_2 \cup \mc M_1$ on $\basis$ using \cref{acyclic_algo} 
	\If{$\mc M_1$ is complete}
		\State set $\block{\bm h}{\mc C_1} = \bm 0$ 
		\State determine $\block{\bm h}{\mc D_1}$ by back substitution using \cref{back_algo}
		\State \textbf{return} $(\block{\bm h}{\mc D_1},\block{\bm h}{\mc C_1})^T$
	\Else
		\State compute new basis $\basis'$ from $\basis$ by change of basis associated with $\mc M_2 \cup \mc M_1$
		\If{$\mc M_1$ is empty}
		\State determine $\block{\bm h}{\mc C_1}$ in \cref{residual} using a sparse linear system solver
		\State \textbf{return} $\block{\bm h}{\mc C_1}$
		\Else
			\State recursively call this routine with input the set of critical elements $\mc C$ of $\basis'=\mc U \cup \mc D \cup \mc C$
			\State set $\block{\bm h}{\mc C_1}$ to the value returned by the recursive process
			\State determine $\block{\bm h}{\mc D_1}$ by back substitution using \cref{back_algo}
		\State \textbf{return} $(\block{\bm h}{\mc D_1},\block{\bm h}{\mc C_1})^T$  
		\EndIf
	\EndIf		
\EndProcedure
\end{algorithmic}
\label{recursive_algo}
\end{algorithm}

For a given input basis $\basis$, \cref{recursive_algo} produces a sequence of bases $\basis=\basis^{\,(0)},\basis^{\,(1)},\dots, \basis^{\,(n)}$, where index $i$ keeps track of the recursion depth.
Each $\basis^{\,(i+1)}$ is constructed from $\basis^{\,(i)}$ by recursion as follows.
We have the partition of $\basis'^{\,(i)} = \basis^{\,(i)} \cdot (\mc M_2^{\,(i)} \cup \mc M_1^{\,(i)}) =  \mc U^{\,(i)} \cup \mc D^{\,(i)} \cup \mc C^{\,(i)}$, associated with acyclic matching $\mc M_2^{\,(i)} \cup \mc M_1^{\,(i)}$ for $i\in \{0,\dots,n\}$.
Note that from the discussion at the end of \cref{greedy} it follows that $\mc M_2^{\,(i)}$ is not empty only for $i=0$.
We set $\basis^{\,(i+1)} = \mc C^{\,(i)}$.

It is clear that the cardinality of each basis $\basis^{\,(i)}$ decreases as long as there are some collapses to be made by \cref{acyclic_algo}.
If no free or flat basis elements are available in \cref{acyclic_algo}, then the recursion is stopped and a sparse linear system solver is used.
So the algorithm always terminates.

\begin{figure}[h]
    \centering
    \includegraphics[scale=0.8]{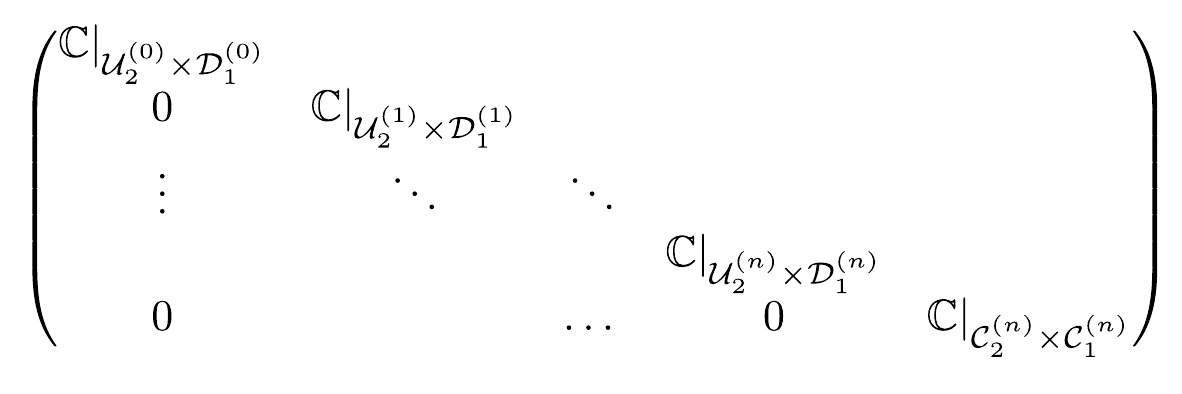}
    \caption{General block structure of matrix $\mat C$ produced by \cref{recursive_algo}.}
    \label{recursive_img}
\end{figure}

At the end of the recursion, we can determine a discrete vector potential $\bm h$ solution of \cref{main} by recursively applying the reasoning described in \cref{residual_sec}.
We write the general form a discrete vector potential $\bm h$ solution of \cref{main} as
\begin{equation}
\bm h = (\block{\bm h}{\mc D_1^{\,(0)}}, \dots, \block{\bm h}{\mc D_1^{\,(n)}}, \block{\bm h}{\mc C_1^{\,(n)}})^T,
\end{equation}
where each $\block{\bm h}{\mc D_1^{\,(i)}}$ is the subvector of $\bm h$ induced by $\mc D_1^{\,(i)}$, see \cref{recursive_img}.
We start by possibly determining $\block{\bm h}{\mc C_1^{\,(n)}}$ using a linear system solver.
Next, we determine each $\block{\bm h}{\mc D_1^{\,(0)}},\dots, \block{\bm h}{\mc D_1^{\,(n)}}$ by back substitution starting from $\block{\bm h}{\mc D_1^{\,(n)}}$ down to $\block{\bm h}{\mc D_1^{\,(1)}}$.
This proves the correctness of \cref{recursive_algo}.

\begin{theorem}
Given as input the canonical basis $\Cbasis$, \cref{recursive_algo} returns a discrete vector potential $\bm h$ solution of \cref{main}.
\end{theorem}

Concerning the computational complexity of \cref{recursive_algo}, if at a certain recursion stage no free or flat basis elements are available in \cref{acyclic_algo}, a linear system solver is employed.
Thus, the worst-case complexity is cubical with respect to the size of the input mesh.
Yet, the average complexity of \cref{recursive_algo} in all tested problems has been linear.
This is because, in practice, only one recursively call of \cref{recursive_algo} is needed.
Indeed, the first recursive call of \cref{recursive_algo} always finds a complete acyclic matching $\mc M_1$ with respect to new basis $\basis^{\,(1)}$.
In other words, there is no need in practice to solve a linear system of the form \cref{residual} with matrix $\block{\mat C}{\mc C_2^{\,(1)} \times \mc C_1^{\,(1)}}$ but instead we can recursively apply back substitution to determine $\block{\bm h}{\mc D_1^{\,(1)}}, \block{\bm h}{\mc D_1^{\,(0)}}$, in this order, after setting free variables in $\mc C_1^{\,(1)}$ to zero.

\section{Numerical results}
\label{numerical}
In this section we illustrate the performance of \cref{recursive_algo}.
We consider different sets of test problems.
In the first set, we focus on simple triangulations that appear in practical boundary value problems.
Next, we present more complicated benchmark triangulations.

For all triangulations that appear in practical boundary value problems, our greedy procedure in \cref{acyclic_algo} always finds a complete acyclic matching $\mc M_1^{\,(0)}$.
Only for the more complicated benchmark problems it is necessary to exploit the novel recursive procedure in \cref{recursive_algo}.
However, to compute a discrete vector potential it is enough, in all tested problems, only one recursively call of \cref{recursive_algo}.

The algorithm has been implemented in C++.
All the numerical computations have been performed in a Intel Core i7-3720QM, with a processor at 2.60 GHz in a laptop with 16 GB of RAM.

\subsection{Triangulations coming from real case boundary value problems}
\label{real_case}
We consider a triangulation coming from a computational electromagnetics application.
As an example, the computation of the source magnetic field for the TEAM problem 7 has been addressed \cite{Fujiwara1990ResultsFB}.

\Cref{real_table} contains information on the number of cells of triangulations of different sizes together with the time (in milliseconds) required to compute the discrete vector potential $\bm h$ using \cref{recursive_algo}.
It is worth noticing that in a triangulation with about 2 million tetrahedra our procedure computes a discrete vector potential under a second.
We run our \cref{acyclic_algo} with different triangulations of the same metal plate and on each example \cref{acyclic_algo} finds a complete acyclic matching $\mc M_1^{\,(0)}$.
Thus, there is no recursive call of \cref{recursive_algo}.
For the considered examples in \cref{real_table} we can clearly see the linear behaviour of the computational time with respect to the size of the triangulations.

\begin{table}[h]
\centering
\caption{Running times of \cref{recursive_algo} for the modified TEAM benchmark example for triangulations of  decreasing size.}
\label{real_table}
\begin{tabular}{ |c||c|c|c|c||c|c| }
 \hline
   Name & Tetrahedra $(\cc)$ & Faces $(\fc)$ & Edges $(\ec)$ & Vertices $(\vc)$ & Time [ms] \\
 \hline
 Mesh 1  & 1,851,493  &   3,871,379  &  2,419,350  & 399,465 &  992 \\
 \hline
 Mesh 2  & 1411688  &   2847256  &  1683787  & 248220 &  756\\
 \hline 
  Mesh 3  & 529,664  &   1,065,104  &  626,566  & 91,127 &  284 \\
 \hline
   Mesh 4  & 186264  &   378588  &  226584  & 34261 &  101 \\
 \hline
\end{tabular}
\end{table}

\subsection{Bing's House}
\label{bing_results}
\begin{figure}[t]
    \centering
    \includegraphics[scale=0.3]{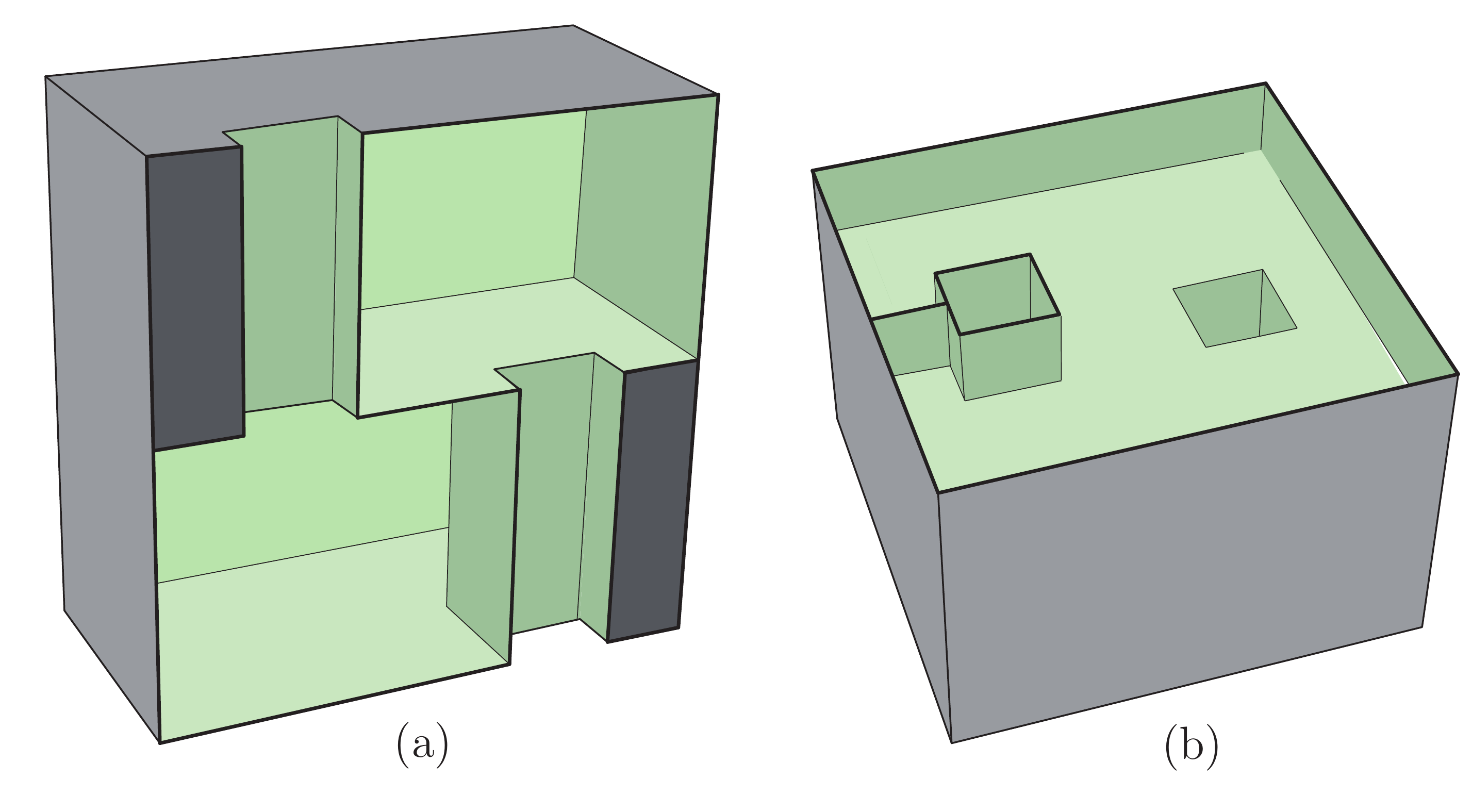}
    \caption{Two views of the considered 3-dimensional thickening of a Bing's house with two rooms.}
    \label{bing_fig}
\end{figure}

A Bing's House is now considered \cite{Cohen1973ACI}.
The simplicial complex, homeomorphic to a 3-dimensional ball, can be obtained by replacing every surface in the Bing's House by a thick wall made of 3-cells.
At the end of this procedure we obtain the polyhedron in \cref{bing_fig}.
Although we can informally identify two ``chambers", it can be demonstrated that the Bing's House is homeomorphic to the three-dimensional ball.

As in the previous set of tests, \Cref{bing_table} contains information about the number of cells of the considered triangulations together with the computational time required to compute a discrete vector potential.
We have found that in almost all runs of \cref{acyclic_algo} we get a complete acyclic matching $\mc M_1^{\,(0)}$.
Only in a few cases we need to resort to a recursive call of \cref{recursive_algo}.
However, to compute a discrete vector potential it is enough, in all these cases, only one recursive call of \cref{recursive_algo}, given that in the first recursion we always find a complete acyclic matching $\mc M_1^{\,(1)}$. 

To measure the complexity of the first recursive call, we consider the cardinality of the basis $\basis_2^{\,(1)}$, namely the output basis of 2-chains becoming the input for the first recursive call of \cref{recursive_algo}.
We have found that the cardinality of  $\basis_2^{\,(1)}$ is always less than 10 on thousands of algorithm runs with different choices of the acyclic matching $\mc M_2^{\, (0)}$ of 2-chains.
Accordingly, as reported in \cref{bing_table}, we observe no influence of the recursive call in \cref{recursive_algo} on the linear behaviour of the running times with respect to the size of the triangulations.

\begin{table}[h]
\centering
\caption{Running times of \cref{recursive_algo} for various triangulations of the  thick Bing's House.}
\label{bing_table}
\begin{tabular}{ |c||c|c|c|c||c|c| }
 \hline
   Name & Tetrahedra $(\cc)$ & Faces $(\fc)$ & Edges $(\ec)$ & Vertices $(\vc)$ & Time [ms] \\
 \hline
 Bing 1  & 800,020  &   1,600,537  &  937,631  & 137,115 &  429 \\
 \hline
  Bing 2  & 87,221  &   175,317  &  102,212 & 14,117 &  47 \\
 \hline
\end{tabular}
\end{table}

\subsection{Knot-theoretic obstructions}
We consider 3-balls of $\R^3$ which admit non-collapsible triangulations.

Obstructions coming from short knots have been considered first in the works \cite{Bing,Goodrick1968NonsimpliciallyCT}.

\begin{figure}[t]
    \centering
    \includegraphics[scale=0.6]{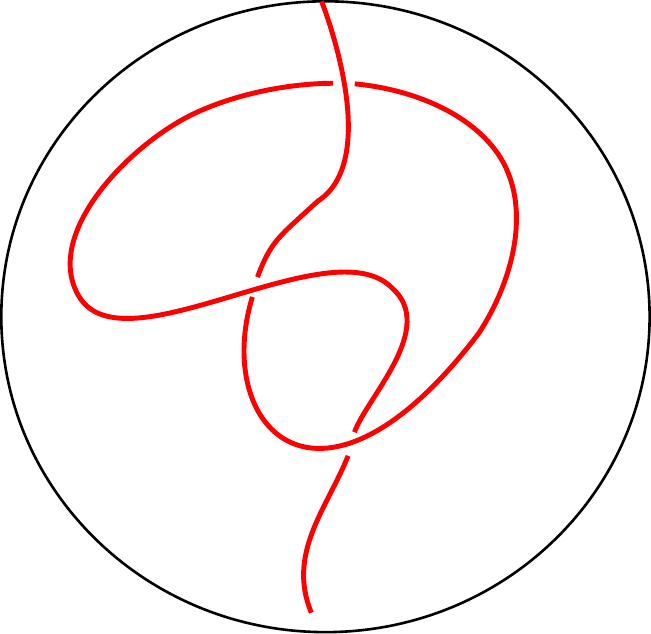}
    \caption{A knotted spanning arc in a 3-ball $\Omega$ at the core of Furch's construction.}
    \label{furch_knot}
\end{figure}

In \cite{Bing}, Bing proved, using knot theory, that some triangulations of the 3-ball are not collapsible.
Bing’s construction works as follows.
One starts with a triangulated 3-ball $\Omega$ and introduces a ``knotted spanning arc" in its 1-skeleton.
A knotted spanning arc is an arc as in \cref{furch_knot}.
We dig a knot-shaped tubular hole inside $\Omega$ starting from the top and we stop digging one step before the tunnel go through the bottom of $\Omega$.
In this way we obtain a 3-ball $\Omega'$ containing a knot having all its edges on the boundary of $\Omega'$, except for a single interior edge.

If the knot is sufficiently complicated (like a double, or a triple trefoil), Bing's ball cannot be collapsible \cite{Bing,Lutz2013KnotsIC}.
In contrast, if the knot is simple enough (like a single trefoil), then Bing's ball may be collapsible.
The construction also appears in the 1924 work \cite{furch} of Furch and for the present discussion we refer to it as \emph{Furch's knotted ball}; see \cite{Ziegler1998ShellingP3} (Section 3.1) for an historical account.

Firstly, we consider the simplest case of Furch's knotted ball with only one trefoil knot.
\Cref{furch_table} summarizes the geometrical information of this triangulation named Furch 1.
In each run of \cref{recursive_algo}, we do not find a complete acyclic matching $\mc M_1^{\,(0)}$ of 1-chains.
Thus, \cref{recursive_algo} is recursively called.
However, only one recursive call is needed in all the considered runs of \cref{recursive_algo}, given that in the first recursion we always find a complete acyclic matching $\mc M_1^{\,(1)}$.

To measure the complexity of the first recursive call, we consider, as done in \cref{bing_results}, the cardinality of the basis $\basis_2^{\,(1)}$.
We have found that the cardinality of  $\basis_2^{\,(1)}$ is always less than 30 on thousands of algorithm runs with different choices of the acyclic matching $\mc M_2^{\, (0)}$ of 2-chains.

As the last benchmark problem, we consider a more complicated obstruction.
We dig one hundred trefoil knots in a parallel-like fashion starting from the top of $\Omega$.
\Cref{furch_table} reports again the results for this triangulation named Furch 100.
We found that the cardinality of  $\basis_2^{\,(1)}$ is always less than 360 on thousands of algorithm runs with different choices of the acyclic matching $\mc M_2^{\, (0)}$ of 2-chains.
Also in this case, only one recursive call is needed in all considered runs of \cref{recursive_algo}, given that in the first recursion we always find a complete acyclic matching $\mc M_1^{\,(1)}$.

We observe no influence of the recursive call in \cref{recursive_algo} on the linear behaviour of the running times  with respect to the size of the triangulations. 
This is because of the small cardinality of $\basis_2^{\,(1)}$.

This example show the effectiveness and generality of our method.
In fact, using the approach in \cite{Rodrguez2015FiniteEP}, we have computed more than one hundred double integral evaluations.
Similarly, using the approach in \cite{Dlotko2011EfficientGS}, we have constructed and then solved a linear system of more than one hundred equations having as unknowns the symbolic variables employed in the approach.
Therefore, Furch 100 is an explicit example of a triangulation on which the approaches \cite{Rodrguez2015FiniteEP,Dlotko2011EfficientGS} perform poorly compared to \cref{recursive_algo}.
To have a provably good method which is reliable in practice, \cref{recursive_algo} is expressly needed.

An important point is that in all tested problems we do not see a dependence between the cardinality of the basis $\basis_2^{\,(1)}$ and the number of cells of the input triangulations.
Thus, this quantity can be used as an indicator of how easy it is to find a complete acyclic matching $\mc M_1^{\,(0)}$ on a given input triangulation, namely, it quantifies the ``topological complexity" of the triangulation.

\begin{table}[h]
\centering
\caption{Running times of \cref{recursive_algo} for the Furch's knotted balls.}
\label{furch_table}
\begin{tabular}{ |c||c|c|c|c||c|c| }
 \hline
   Name & Tetrahedra $(\cc)$ & Faces $(\fc)$ & Edges $(\ec)$ & Vertices $(\vc)$ & Time [ms] \\
 \hline
  Furch 100  & 243062  &   506619  &  311547 & 47991 &  129 \\
 \hline
 Furch 1  &  31189 &   63830  &  38445  & 5805 &  17 \\
 \hline
\end{tabular}
\end{table}

\section{Conclusions}
\label{conclusions}
The novel algorithm presented in this paper was proved to be general, straightforward to implement and outperforms competing state-of-art algorithms in the class of admissible meshes while showing linear average complexity with respect to the input mesh size.
By applying discrete Morse theory, we have shown that for the important class of
simplicial triangulations we achieve linear computational complexity for all considered test problems.
These include, besides real case triangulations having simple topological properties, also pathological triangulations.
A challenging test case made of one hundred trefoil knots have been considered and yet the proposed algorithm succeeds in computing a discrete vector potential.
Thus, we expect that our algorithm works for every practical mesh that one encounters in practical applications.

Worst-case complexity analysis can be misleading in the analysis of our \cref{recursive_algo}.
Worst-case complexity analysis of our algorithm assumes that at certain iteration no new matched pairs are possible.
In this case we need to employ a sparse linear system solver, which leads to a cubical worst-case complexity with respect to mesh size parameters.
Yet, the average complexity is linear in all tested problems.
Our point of view is that the reason we see  linear computational complexity is because our examples restrict to 3-dimensional cell complexes decomposing bounded domains of $\R^3$ with sufficiently regular boundaries.
We have observed that our recursive strategy based on algebraic discrete Morse theory is efficient at solving potential topological obstructions that can appear in 3-dimensional cases.
However, to prove that linear worst-case times are guaranteed is still an open problem and will be subject to future work.

\bibliographystyle{model1-num-names}
\bibliography{jcp}

\end{document}